\UseRawInputEncoding
\documentclass[12pt]{article}

\usepackage[T1]{fontenc}        
\usepackage[utf8]{inputenc}     
\usepackage{lmodern}            
\usepackage{amsmath, amssymb, amsthm}

\usepackage{amsthm}

\newtheorem{theorem}{Theorem}
\newtheorem{proposition}{Proposition}

\usepackage{pgfplotstable}      
\usepackage{booktabs}           
\usepackage{array}
\usepackage{longtable}

\usepackage{graphicx}
\usepackage{float}

\usepackage{geometry}
\usepackage{listings}
\usepackage{xcolor}

\usepackage{hyperref}
\hypersetup{
  pdfborder={0 0 0},
  colorlinks=true,
  linkcolor=black,
  citecolor=black,
  urlcolor=black,
  breaklinks=true
}

\newcommand{\Erdos}{Erd\H{o}s}

\title{Partial Resolution of the Erdős--Straus, Sierpiński, and Generalized Erdős--Straus Conjectures Using New Analytical Formulas}

\author{
  Philemon Urbain MBALLA\\[0.5em]
  \texttt{philemonmballa@gmail.com}\\
  \texttt{philemon-urbain.mballa@etu.u-paris.fr}
}

\date{}
\pgfplotsset{compat=1.18}
\begin{document}

\maketitle

\begin{abstract}
This article proposes a unified analytical approach leading to a partial resolution of the Erdős–Straus, Sierpiński conjectures, and their generalization. We introduce an equivalent reformulation of these conjectures while constructing two new explicit analytical formulas. The first formula, which is a special case of the second, is based on a divisibility condition, whereas the second, more general formula, relies on the existence of a perfect square, which we conjecture to always hold. Under these conditions, the formulas verify the conjectures even for very large numerical values. Moreover, our method reduces the problem to the search for a suitable perfect square, thereby opening the way to a complete proof of these conjectures. In conclusion, we present open questions and conjectures to the mathematical community regarding the generalization of these formulas.
\end{abstract}

\section{Introduction}

The Erdős--Straus conjecture, formulated in 1948 by Paul Erdős and Ernst G. Strauss, states that for every natural number $n \geq 2$, there exist $x, y, z \in \mathbb{N}^*$ such that
\begin{equation}
    \frac{4}{n} = \frac{1}{x} + \frac{1}{y} + \frac{1}{z}.
\end{equation}

Louis Mordell demonstrated that if $n$ is not congruent to $1, 11^2, 13^2, 17^2, 19^2,$ or $23^2$ modulo $840$, then the conjecture holds. Furthermore, numerical verifications have confirmed its validity up to $n = 10^{17}$.

For decades, mathematicians have been intrigued by the following question:
\textbf{Does an explicit formula exist relating $x, y, z$, and $n$ that systematically verifies the conjecture?}

In this article, we propose two new analytical formulas obtained through an in-depth study of the Diophantine equation~(1). The first relies on a natural divisibility condition, while the second is based on the existence of a perfect square. These results provide a partial proof of the conjecture and represent a significant advancement in its study.

The Polish mathematician Wacław Sierpiński also investigated this conjecture and posed a similar question:
\textbf{Is a similar decomposition possible if $4$ is replaced by $5$?}
He conjectured that for every natural number $n \geq 2$, there exist $x, y, z \in \mathbb{N}^*$ such that
\begin{equation}
    \frac{5}{n} = \frac{1}{x} + \frac{1}{y} + \frac{1}{z}.
\end{equation}

In this article, we will demonstrate that our formulas also verify this conjecture under the same conditions as before.

Finally, following the formulation of the Erdős--Straus conjecture, several mathematicians explored the possibility of generalizing this Egyptian fraction representation. This led to the generalized Erdős--Straus conjecture, which may be more precisely formulated as follows: for every integer $a \geq 5$, there exists an integer $N_0 \in \mathbb{N}$ with $N_0 \geq 2$ (possibly depending on $a$, since it is easy to see that  $N_0$ cannot be strictly smaller than $a/3$)
 such that for all $n \geq N_0$, there exist $x, y, z \in \mathbb{N}^*$ satisfying
\begin{equation}
    \frac{a}{n} = \frac{1}{x} + \frac{1}{y} + \frac{1}{z}.
\end{equation}

Moreover, this result is obtained using the same analytical formulas developed for the Erdős--Straus and Sierpiński conjectures.

In this article, we present a detailed analytical approach leading to these formulas, along with numerical tests supporting our results.

\section{The Erdős-Straus Conjecture}

One of the methods that is known nowaday is that if $n$ is even integer , meaning $n = 2k$ with $k \in \mathbb{N}^*$. In this case, it is well established that:
\begin{equation}
\frac{4}{n} = \frac{1}{\frac{n}{2}} + \frac{1}{n} + \frac{1}{n}. \tag{4}
\end{equation}

This result is one of the most trivial among those obtained since the conjecture was stated, as it shows that all even numbers satisfy the conjecture.

In our new approach, I will demonstrate how this result can be naturally recovered when $n$ is even. Before introducing the general method that allowed me to establish two formulas verifying the conjecture under a certain condition, I will first illustrate my process by analytically testing the conjecture for some values of $n$.

 Our method is based on the structural analysis of equation (1). It is a simple yet powerful and general approach, applicable to even numbers, odd numbers, and prime numbers alike.

---

\subsection{\texorpdfstring{Case $n = 2$ (Even Number)}{Case n = 2 (Even Number)}}

The fundamental equation becomes:
\[
\frac{4}{2} = \frac{1}{x} + \frac{1}{y} + \frac{1}{z}.
\]
I choose not to simplify $\frac{4}{2}$ immediately in order to maintain the general structure of the method. However, **one can simplify this fraction directly**, and by following the same steps that I will present, the same result is always obtained.

The equation can then be rewritten as:
\[
\frac{4}{2} - \frac{1}{x} = \frac{1}{y} + \frac{1}{z},
\]
which reformulates as:
\[
\frac{4x - 2}{2x} = \frac{y+z}{yz}.
\]
Since we are searching for solutions in $\mathbb{N}^*$, it is necessary that $4x - 2$ be strictly positive, that is, $4x - 2 > 0$, which gives $x \geq \lfloor \frac{1}{2} \rfloor + 1 = 1$.

For the rest of the document, **the notation $\lfloor x \rfloor$ will always denote the integer part of $x$** to avoid any confusion.

We then obtain the equation:
\[
 yz(4x - 2) = 2x(y + z).
\]
At this stage, I apply a method commonly used in geometry of curves and surfaces, called parameterization.

---

\subsection{Parameterization Technique}

Let $a,  c \in \mathbb{R}$,and $b,d\in\mathbb{R^*}$ such that:
\[
 ab = cd.
\]
Then, for  $t \in \mathbb{R}$, we can set:
\[
 a = d t, \quad c = b t.
\]
This parametric representation holds for  $t \in \mathbb{R}$ and allows rewriting $(d t) b = (b t) d$, which verifies the initial equation. This technique will be crucial for my results.

---

\subsection{Application to Our Equation}

We then set:
\[
 y + z = (4x - 2) t, \quad yz = 2x t.
\]
This parameterization is valid for  $t \in \mathbb Q{^*}_{+}$.

According to the sum and product method, **$y$ and $z$ are solutions of the quadratic equation:**
\[
 V^2 - ((4x - 2) t) V + 2x t = 0.
\]
Its discriminant is given by:
\[
 \Delta = t^2 (4x - 2)^2 - 8x t.
\]
We imposed $x \geq 1$. Choosing $x = 1$, we obtain:
\[
 \Delta = 4t^2 - 8t = t(4t - 8).
\]
For $y$ and $z$ to be integers, $\Delta$ must be a perfect square, and $4x - 2$ must have the same parity as $\sqrt{\Delta}$.

Taking $t = 2$ (a trivial value for which $\Delta = 0$), we obtain a double root:
\[
 y = z = \frac{(4 \times 1 - 2) \times 2}{2} = 2.
\]
Thus, for $n = 2$, we find $x = 1, y = 2, z = 2$, verifying:
\[
 \frac{4}{2} = \frac{1}{1} + \frac{1}{2} + \frac{1}{2}.
\]
This confirms the validity of the conjecture for $n = 2$.

---

\subsection{Analysis of the Method}

if  we had directly applied the known formula for even numbers given by equation  ($4$) , we would have found:
\[
 x = \frac{2}{2} = 1, \quad y = z = 2.
\]
Thus, we obtain exactly the same values.

One of the key advantages of our analytical approach is its applicability not only to even numbers but also to odd numbers and prime numbers. Unlike classical methods—which have relied primarily on brute-force techniques that merely test the existence of solutions without offering insight into the problem’s structure—our method derives solutions directly from the inherent logical structure of the underlying Diophantine equation. This provides a general, rigorous framework for addressing the Erdős–Straus, Sierpiński, and generalized Erdős–Straus conjectures. Moreover, the elegance and clarity of our approach make it a valuable pedagogical tool for teaching methods of solving Diophantine equations, as it not only yields the solutions but also deepens our understanding of the arithmetic properties involved

Of course, testing the conjecture for specific values of $n$ does not constitute a general proof. My goal here is to demonstrate the efficiency of this method in solving equation (1) analytically.

Next, we will examine the cases $n = 3$, $n = 17$, and finally $n = 841$ (one of the numbers that Mordell had not proven), before moving on to generalization. Since the main challenge concerns odd numbers, I will not dwell much on even numbers. However, you can reproduce this reasoning for any even $n$, and it will work the same way.

\subsection{\texorpdfstring{Case $n = 3$ (odd and prime Number)}
{Case n = 3 (odd and prime number)}}

Applying the same procedure, we impose:
\[
 x \geq \left\lfloor \frac{3}{4} \right\rfloor +1 = 1.
\]
We then obtain the system:
\[
\begin{cases}  
y + z = (4x - 3)t, \\  
yz = 3xt,  
\end{cases}
\]
for  $t \in \mathbb Q{^*}_{+}$.  

The values of $y$ and $z$ satisfy the quadratic equation:
\[
V^2 - ((4x - 3)t)V + 3xt = 0.
\]
The discriminant is given by:
\[
\Delta = t^2 (4x - 3)^2 - 12xt.
\]
We choose $x \geq \left\lfloor \frac{4}{3} \right\rfloor +1 = 2$.

Taking $x = 2$, we obtain:
\[
\Delta = 25t^2 - 24t = t(25t - 24).
\]
A trivial solution ensuring that $\Delta$ is a perfect square is $t = 1$, yielding $\Delta = 1$.

Thus, we find:
\[
y = \frac{(4 \times 2 - 3) - 1}{2} = 2, \quad z = \frac{(4 \times 2 - 3) + 1}{2} = 3.
\]
For $n = 3$, we have $x = 2$, $y = 2$, and $z = 3$, verifying:
\[
\frac{4}{3} = \frac{1}{2} + \frac{1}{2} + \frac{1}{3}.
\]
By selecting a different $x$, we can determine corresponding values of $t$ that yield integer solutions. 
For instance, choosing $x = 4$ results in:
\[
\Delta = 169t^2 - 48t = t(169t - 48).
\]
For $t = 1$, we find a perfect square $11^2 = 121$, yielding $y = 1$ and $z = 12$.

---

\subsection{\texorpdfstring{Case $n = 17$ (odd and prime Number)}
{Case n = 17 (odd and prime number)}}

We impose:
\[
 x \geq \left\lfloor \frac{17}{4} \right\rfloor +1 = 5.
\]
The system obtained is:
\[
\begin{cases}  
y + z = (4x - 17)t, \\  
yz = 17xt,  
\end{cases}
\]
for  $t \in \mathbb Q{^*}_{+}$.  

$y$ and $z$ satisfy the quadratic equation:
\[
V^2 - ((4x - 17)t)V + 17xt = 0.
\]
The discriminant is:
\[
\Delta = t^2(4x - 17)^2 - 68xt.
\]
To minimize $(4x - 17)^2$ and simplify calculations, we choose the smallest $x \geq 5$, i.e., $x = 5$.

This gives:
\[
\Delta = 9t^2 - 340t = t(9t - 340).
\]
We require $\Delta$ to be a perfect square, meaning $9t - 340 \geq 0$, so $t \geq \frac{340}{9} \approx 37.77$.

Using a Python script to scan $t$ in the range $[38, 1000]$, we find the perfect square values of $\Delta$:
\[
 t = 68 \Rightarrow \sqrt{\Delta} = 136, \quad t = 180 \Rightarrow \sqrt{\Delta} = 480.
\]
For $t = 68$, we obtain:
\[
 y = 34, \quad z = 170.
\]
Thus, we confirm:
\[
\frac{4}{17} = \frac{1}{5} + \frac{1}{34} + \frac{1}{170}.
\]
The reader may compute $y$ and $z$ for $t = 180$ using the same approach.

---
\subsection{\texorpdfstring{Case $n = 841$ (number not proven by Mordell)}
{Case n = 841 ( Number not proven by Mordell )}}

Following the same procedure, we impose:
\[
 x \geq \left\lfloor \frac{841}{4} \right\rfloor +1 = 211.
\]
The system obtained is:
\[
\begin{cases}  
y + z = (4x - 841)t, \\  
yz = 841xt,  
\end{cases}
\]
for all $t \in \mathbb Q{^*}_{+}$.  

$y$ and $z$ satisfy the quadratic equation:
\[
V^2 - ((4x - 841)t)V + 841xt = 0.
\]
The discriminant is:
\[
\Delta = t^2(4x - 841)^2 - 3364xt.
\]
Choosing $x = 211$, we obtain:
\[
\Delta = 9t^2 - 709804t = t(9t - 709804).
\]
For $\Delta$ to be a perfect square, we require $9t - 709804 \geq 0$, so $t \geq \frac{709804}{9} \approx 78867.11$.

Scanning for $t$ in the range $[78868, 500000]$, we find:
\[
 t = 185600, \quad \sqrt{\Delta} = 422240.
\]
Thus, we compute:
\[
 y = 67280, \quad z = 489520.
\]
Verifying:
\[
\frac{4}{841} = \frac{1}{211} + \frac{1}{67280} + \frac{1}{489520}.
\]

\vspace{0.5em} 
Both sides yield the same numerical value.

\vspace{0.5em}

\noindent\textbf{NB.} In most cases, we restrict the search for $t$ to $\mathbb{N}^*$. However, numerical simulations indicate that allowing $t$ to range over non-integer positive rationals $\mathbb{Q}_+^* \setminus \mathbb{N}$ becomes necessary for $a \geq 5$ when $n$ is small, namely in a neighborhood of the threshold $N_0$, since solutions remain sparse in this region. Beyond a certain threshold, we observe an abundance of integer values of $t$ for which the discriminant is a perfect square. Moreover, $t$ may also naturally take values in $\mathbb{Q}_+^*$ in the context of Formula~(1) as soon as the divisibility condition is satisfied.

---

\section{Conclusion}

This analytical approach effectively solves equation (1) and finds solutions satisfying the conjecture for even, odd, and prime numbers.

With this method, one can manually prove the conjecture for many numbers not previously proven by Mordell.

However, verifying the conjecture for specific values does not constitute a general proof. We will now proceed to the generalization and present the results obtained.

\section{Generalization of the Erdős-Straus Conjecture}

We want to prove that for all \( n \) in \( \mathbb{N}^* \) with \( n \geq 2 \), there exist \( x, y, z \) in \( \mathbb{N}^* \) such that  

\[
\frac{4}{n} = \frac{1}{x} + \frac{1}{y} + \frac{1}{z}
\]

We have called this equation (1).  

My proposed proof follows the same structure as before, as you may have noticed, the value of \( n \) does not change the structure of the proof. I will use exactly the same procedure as previously.  

Equation (1) is equivalent to:  

\[
\frac{4}{n} - \frac{1}{x} = \frac{y + z}{yz}
\]

which is also equivalent to  

\[
\frac{4x - n}{nx} = \frac{y + z}{yz}
\]

For the right-hand side to be strictly positive, it is necessary that \( 4x - n > 0 \), i.e., \( x > \frac{n}{4} \), so we have \( x \geq \left\lceil \frac{n}{4} \right\rceil \).  

We still have:  

\[
\frac{4x - n}{nx} = \frac{y + z}{yz}
\]

which holds if and only if  

\[
yz(4x - n) = nx(y + z)
\]

Similarly, using parameterization, we can express:  

\[
\begin{cases}
y + z = (4x - n)(2t) \\
yz = nx(2t)
\end{cases}
\]

for  $t \in \mathbb Q{^*}_{+}$. If we replace \( y + z \) and \( yz \) with these expressions, the equation \( yz(4x - n) = nx(y + z) \) will always hold.  

Later, I will explain why I parameterized using \( 2t \) in this general case instead of using \( t \) as I did for specific values of \( n \). (Of course, we are allowed to parameterize using any multiple of $t \in \mathbb Q{^*}_{+}$, as the equation will always hold.)  

Now, let us use sums and products. From the system above, \( y \) and \( z \) satisfy the quadratic equation in \( V \):  

\[
V^2 - (\text{sum})V + \text{product} = 0
\]

which gives:  

\[
V^2 - ((4x - n)(2t))V + 2ntx = 0.
\]

The discriminant is given by:  

\[
\Delta = 4t^2(4x - n)^2 - 8ntx. \quad (5)
\]

You may recall that I promised to retrieve equation (4) for even values of \( n \). If we set \( n = 2x \), meaning \( n \) is even since \( x \) is in \( \mathbb{N}^* \), the discriminant becomes:  

\[
\Delta = 4t^2(4x - 2x)^2 - 16x^2t.
\]

\[
= 16x^2t^2 - 16x^2t.
\]

An obvious choice for \( t \) so that \( \Delta \) is a perfect square is \( t = 1 \), which gives \( \Delta = 0 \).  

Thus,  

\[
y = z = \frac{2(4x - 2x)}{2} = 2x
\]

so \( y = z = n = 2x \), which implies \( x = \frac{n}{2} \), and we obtain:  

\[
\frac{4}{n} = \frac{1}{\frac{n}{2}} + \frac{1}{n} + \frac{1}{n}
\]

which retrieves equation (4).  

This was a side note, as our goal is to generalize. So, we return to our quadratic equation in \( V \) in its general form, with \( n \in \mathbb{N}^* \) such that \( n \geq 2 \) (as required by the conjecture).  

We still have our discriminant computed in (5):

To obtain integer solutions for this equation, it is necessary that \( \Delta \) be a perfect square. I struggled to find an obvious perfect square for \( \Delta \). It seems to me that the polynomial given in (5) is rich, but unfortunately, after some quick and unsuccessful attempts to find a classic perfect square of \( \Delta \) without imposing any condition on \( n \) (since we are in the process of generalization), the idea I had was to use the fact that the discriminant \( \Delta \) holds for  $t \in \mathbb{Q}_+^*$ arbitrary and see what can be deduced to obtain a perfect square in order to bypass the difficulty. Let's see how I proceed
Let’s revisit the compact form of \( \Delta \) in (5):
\section{\texorpdfstring{Compact Form of \( \Delta \)}{Compact Form of Delta}}

\[
\Delta = 4t^2(4x - n)^2 - 8ntx
\]
This simplifies to:

\[
\Delta = 4t(t(4x - n)^2 - 2nx).
\]
This is valid for  $t \in \mathbb{Q}_+^*$ arbitrary because, as mentioned earlier, the parametrization is valid for  $t \in \mathbb{Q}_+^*$ arbitrary. You also observed how this information allowed us to prove the conjecture (manually for some \( n \)) by providing specific values of \( t \) that help us obtain a perfect square. We will use this information in this generalization.

We saw that \( \Delta >0 \) for
\[
t > \frac{2nx}{(4x - n)^2}.
\]

 Recalling the expression for \( \Delta \),
\[
\Delta = 4t(t(4x - n)^2 - 2nx),
\]
and since \( \Delta > 0 \)  for \[
 t > \frac{2nx}{(4x - n)^2}.
\]  for the x and t which respect this condition, it is clear that there exist \( x \) and \( t \) such that \( t(4x - n)^2 - 2nx \) is in \( \mathbb{{Q}_+^*}  \). We can then write:
\[
t(4x - n)^2 - 2nx = k,
\]
with $ k \in \mathbb{Q}_+^*$ .
Thus,
\[
\Delta = 4t(t(4x - n)^2 - 2nx) = 4tk,
\]
which holds for  $t \in \mathbb{Q}_+^*$ arbitrary. So, for \( t = k \) (This is why we excluded the case $\Delta = 0$, otherwise $k = \left( t(4x - n)^2 - 2nx \right)$ can be zero  
and we cannot set $t = k$ because from the start, during the parameterization, $t$ is non-zero).
, so we have:
\[
\Delta = 4tk = 4k^2 = (2k)^2.
\]
We have thus obtained a perfect square (rational or integer perfect square). Let’s now see under what conditions this transformation holds.

\section{Conditions for Obtaining a Perfect Square}
To obtain this perfect square, we wrote \( t(4x - n)^2 - 2nx = k \), and then set \( t = k \), simplifying the equation to:
\[
k(4x - n)^2 - 2nx = k,
\]
or:
\[
k\left((4x - n)^2 - 1\right) = 2nx.
\]
This implies that
\[
k = \frac{2nx}{(4x - n - 1)(4x - n + 1)},
\qquad k \in \mathbb{Q}_+^*.
\]
At this stage, no additional condition is required, since the above expression naturally defines a positive rational value of $k$ without any restriction.

Since, \( \Delta = 4k^2 = (2k)^2 \), so:
\[
V = \frac{(4x - n)(2k) \pm 2k}{2}.
\]
there is a simplification of 2, This is why I had previously parameterized by \( 2t \), and we just replaced \( t \) with \( k \) to calculate the roots of this polynomial, and this was possible because we set \( t = k \) to obtain a perfect square.

we have \( y = k(4x - n - 1) \) and \( z = k(4x - n + 1) \). Now, from earlier, we know that:

\[
k = \frac{2n x}{(4x - n - 1)(4x - n + 1)}
\]

Substituting \( k \) into the expressions for \( y \) and \( z \), we get:

\[
y = \frac{2n x}{(4x - n - 1)(4x - n + 1)} (4x - n - 1) = \frac{2n x}{4x - n + 1}
\]

and

\[
z = \frac{2n x}{(4x - n - 1)(4x - n + 1)} (4x - n + 1) = \frac{2n x}{4x - n - 1}
\]

A natural condition for $y$ and $z$ to be integers is that both $4x - n - 1$ and $4x - n + 1$ divide $2nx$.

Next, we find the conditions for \( x \) such that both \( y \) and \( z \) are positives. Since \( 2n x > 0 \), we only need to ensure that the denominators are strictly positive.

\textbf{Case for \( y \):}

\( y > 0 \) if:

\[
4x - n + 1 > 0 \quad \text{which implies} \quad x > \frac{n - 1}{4}
\]

Thus, for \( y \), we need \( x \geq \left\lceil \frac{n - 1}{4} \right\rceil + 1 \).

\textbf{Case for \( z \):}

\( z > 0 \) if:

\[
4x - n - 1 > 0 \quad \text{which implies} \quad x > \frac{n + 1}{4}
\]

Thus, for \( z \), we need \( x \geq \left\lceil \frac{n + 1}{4} \right\rceil + 1 \).

Therefore, the overall condition is:

\[
x \geq \max\left\{\left\lceil \frac{n - 1}{4} \right\rceil + 1, \left\lceil \frac{n + 1}{4} \right\rceil + 1 \right\}
\]

which gives:

\[
x \geq \left\lceil \frac{n + 1}{4} \right\rceil + 1
\]

Thus, \( x \in \left[\left\lceil \frac{n + 1}{4} \right\rceil + 1, +\infty \right[ \).

Now, let us check for which values of \( x \), \( y \), and \( z \) are in \( \mathbb{N}^* \). We recall the following property: if \( a \in \mathbb{N}^* \) and \( b \in \mathbb{N} \), and if \( a \) divides \( b \), then \( a \leq b \).

Therefore, \( y \in \mathbb{N}^* \) if and only if \( 4x - n + 1 \) divides \( 2n x \), which implies:

\[
4x - n + 1 \leq 2n x
\]

\[
4x - 2n x \leq n - 1
\]

\[
x(4 - 2n) \leq n - 1
\]

For \( n > 2 \), \( 4 - 2n < 0 \), so:

\[
x \geq \frac{n - 1}{2n - 4}, \quad x \in \mathbb{N}
\]

Similarly, we can show that \( z \in \mathbb{N}^* \) if:

\[
x \geq \frac{n + 1}{2n - 4}
\]

Thus, we have the following conditions:

\[
x \geq \left\lceil \frac{n + 1}{4} \right\rceil + 1
\]

\[
x \geq \frac{n - 1}{2n - 4}
\]

\[
x \geq \frac{n + 1}{2n - 4}
\]

The two latter values of \( x \) are not valid for \( n = 2 \). Therefore, the correct choice that satisfies all conditions is:

\[
x \geq \left\lceil \frac{n + 1}{4} \right\rceil + 1
\]

In other words, \( x \in \left[\left\lceil \frac{n + 1}{4} \right\rceil + 1, +\infty \right[ \).

As we will see in the numerical checks, under these divisibility conditions, this formula will allow us to decompose the fraction \( \frac{4}{n} \) into three unit fractions for those \( n \) satisfying these conditions up to infinity (since we have shown that the values of \( x \) for which \( 4x - n - 1 \) and \( 4x - n + 1 \) divide \( 2n x \) extend to infinity). This is good news because we always hope that for any given \( n \), there is an infinite number of choices for \( x \) that satisfy the divisibility condition, allowing us to prove the conjecture for that specific \( n \).

Before proceeding with this, let's first verify whether, under these divisibility conditions, the Erdős–Strauss conjecture is satisfied for:

\[
y = \frac{2n x}{4x - n + 1} \quad \text{and} \quad z = \frac{2n x}{4x - n - 1}
\]

We need to check if:

\[
\frac{1}{x} + \frac{1}{y} + \frac{1}{z} = \frac{4}{n}
\]

Substituting the expressions for \( y \) and \( z \):

\[
\frac{1}{x} + \frac{1}{\frac{2n x}{4x - n + 1}} + \frac{1}{\frac{2n x}{4x - n - 1}} = \frac{1}{x} + \frac{4x - n + 1}{2n x} + \frac{4x - n - 1}{2n x}
\]

Simplifying:

\[
\frac{1}{x} + \frac{(4x - n + 1) + (4x - n - 1)}{2n x} = \frac{1}{x} + \frac{8x - 2n}{2n x} = \frac{1}{x} + \frac{4x - n}{n x} = \frac{4}{n}
\]

Thus, \( x \), \( y \), and \( z \) satisfy the Erdős–Strauss conjecture under the established divisibility conditions.

In summary, we have proven that for all \( n \geq 2 \), If there exists \( x_0 \in \mathbb{N}^* \) with \( x_0 \in \left[ \left\lceil \frac{n+1}{4} \right\rceil + 1, +\infty \right[ \).
 such that \( 4x_0 - n - 1 \) and \( 4x_0 - n + 1 \) divide \( 2n x_0 \), then there exist \( y = \frac{2n x_0}{4x_0 - n - 1} \) and \( z = \frac{2n x_0}{4x_0 - n + 1} \) in \( \mathbb{N}^* \) such that:

\[
\frac{4}{n} = \frac{1}{x_0} + \frac{1}{y} + \frac{1}{z}
\]
we therefore have the following theorem:

\begin{theorem}

Let $n \geq 2$ be an integer. Suppose there exists $x_0 \in \mathbb{N}^*$ with
\[
x_0 \in \left[ \left\lceil \frac{n+1}{4} \right\rceil + 1, +\infty \right)
\]
such that both integers $4x_0 - n - 1$ and $4x_0 - n + 1$ divide $2n x_0$.
Then there exist $y, z \in \mathbb{N}^*$ given by
\[
y = \frac{2n x_0}{4x_0 - n - 1}
\qquad \text{and} \qquad
z = \frac{2n x_0}{4x_0 - n + 1},
\]
such that
\[
\frac{4}{n} = \frac{1}{x_0} + \frac{1}{y} + \frac{1}{z}.
\]
\end{theorem}

Below, I will provide numerical calculations that, depending on the range searched for \( x_0 \), cover a some $n$  . Based on this formula verifying the conjecture, you will also be able to perform calculations using more powerful computational resources to identify very large \( x_0 \) and \( n \). Of course, with numerical computations, it is impossible to reach \( +\infty \), but what assures us of the validity of the results is the analytical result we have just established. Every time we find such an \( x_0 \), the conjecture holds for the corresponding \( n \).

My question to the mathematical community for a complete proof of the conjecture is as follows: How can we explicitly find \( x_0 \)? If it is not possible to obtain an explicit expression for \( x_0 \) (which may depend on \( n \)), can we prove the existence of \( x_0 \) by induction on \( n \) (whether classical or strong)? I remain available for any suggestions or collaborations to address this question and definitively resolve the conjecture.

\section{\textbf{Second Formula: Condition of the Existence of a Perfect Square}}

While the first formula is derived from a divisibility condition, the second one is more fundamental in nature, as it provides a genuinely more general framework. Indeed, it is based on the existence of a perfect square and subsumes the first formula as a particular case. We now introduce this second formula.

Let us reconsider the discriminant \( \Delta \) of our quadratic equation:
\[
\Delta = 4t^2(4x - n)^2 - 8ntx = 4\bigl(t^2(4x - n)^2 - 2ntx\bigr).
\]
We have already seen that $\Delta \geq 0$ for $t > \frac{2nx}{(4x - n)^2}$ with $t \in \mathbb{N}^*$. Therefore, for such values of $t$, we define

$t^2(4x - n)^2 - 2ntx = q^2$, with $q \in \mathbb{N}$ (as mentioned previously, here and in most cases we take $t$ and $q$ to be integers in order to automatically guarantee that, once the perfect square condition is satisfied, the resulting solutions are integral. Our parametrization allows $t$ to take arbitrary values without affecting the structure of the equation, and it is therefore up to us to choose suitable values of $t$ ensuring that the corresponding values of $y$ and $z$ are integers).

Thus, we have:
\[
\Delta = 4q^2 = (2q)^2.
\]
We will now prove that the conjecture is true under the condition that this perfect square exists.  

Let \( y = t(4x - n) - q \) and \( z = t(4x - n) + q \). We need to show that, with this formula, the conjecture holds for those \( q \) such that  \[
t^2(4x - n)^2 - 2ntx = q^2, \quad q \in \mathbb{N}.
\].

It is sufficient to check if
\[
\frac{1}{x} + \frac{1}{y} + \frac{1}{z} = \frac{1}{x} + \frac{1}{t(4x - n) - q} + \frac{1}{t(4x - n) + q}.
\]
We can simplify this expression:
\[
\frac{1}{x} + \frac{1}{t(4x - n) - q} + \frac{1}{t(4x - n) + q} = \frac{1}{x} + \frac{2t(4x - n)}{t^2(4x - n)^2 - q^2}.
\]
Now, since \( q^2 = t^2(4x - n)^2 - 2ntx \), we have:
\[
\frac{1}{x} + \frac{2t(4x - n)}{t^2(4x - n)^2 - t^2(4x - n)^2 + 2ntx} = \frac{1}{x} + \frac{2t(4x - n)}{2nxt}.
\]
Simplifying further:
\[
= \frac{1}{x} + \frac{4x - n}{nx} = \frac{n - n + 4x}{nx} = \frac{4}{n}.
\]
Thus, the second formula allows us to state the following theorem.

\begin{theorem}
For any $n \in \mathbb{N}$ with $n \geq 2$, if there exists $(x_0,t_0) \in (\mathbb{N}^*)^2$ such that
\[
t_0^2(4x_0 - n)^2 - 2nt_0x_0 = q^2 \in \mathbb{N},
\]
then there exist
\[
y = t_0(4x_0 - n) - q \quad \text{and} \quad z = t_0(4x_0 - n) + q,
\]
with $(y,z) \in \mathbb{N}^2$, such that
\[
\frac{4}{n} = \frac{1}{x_0} + \frac{1}{y} + \frac{1}{z}.
\]
\end{theorem}

This also proves the conjecture under this condition. I will also present below the numerical values of \( n \) found within a range where we search for \( (x_0, t_0) \) through numerical calculations.

I also pose the following question to the mathematical community: how can we find an explicit formula for \( (x_0, t_0) \) that works for all \( n \) without any conditions? In other words, how can we prove that the polynomial \( t^2(4x - n)^2 - 2ntx \) or \( 4t^2(4x - n)^2 - 8ntx \) always admits a perfect square for all \( n \in \mathbb{N} \), with \( n \geq 2 \), in order to complete the proof of the conjecture?

I remain available for suggestions or collaborations to address this question.

This is the contribution of this article to the Erdős-Strauss conjecture.

\section{ Sierpi\'nski Conjecture}

The Sierpi\'nski conjecture is similar to the Erd\H{o}s-Straus conjecture, except that 4 is replaced by 5. It states that for any $n \in \mathbb{N}$ with $n \geq 2$, there exist $x, y, z \in \mathbb{N}^*$ such that:

\begin{equation} \label{eq:seven}
    \frac{5}{n} = \frac{1}{x} + \frac{1}{y} + \frac{1}{z}.\tag{7}
\end{equation}

A trivial case arises when $n$ is divisible by 3. In this case, one can simply take $x = \frac{n}{3}$, $y = z = n$, and the conjecture is verified for all $n$ divisible by 3.

The good news is that with this new approach, everything stated for the Erd\H{o}s-Straus conjecture remains valid for the Sierpi\'nski conjecture by simply replacing 4 with 5. This demonstrates the power and generality of this method.

To keep the article concise, I will not manually test specific values of $n$ as I did for the Erd\H{o}s-Straus conjecture. However, the reasoning remains the same, and the reader may verify this independently.

\subsection{Transforming the Equation}

Equation (7) is equivalent to:

\begin{equation}
    \frac{5}{n} - \frac{1}{x} = \frac{y+z}{yz},
\end{equation}

which can be rewritten as:

\begin{equation}
    \frac{5x - n}{nx} = \frac{y+z}{yz}.
\end{equation}

For the right-hand side to be strictly positive, it is necessary that $5x - n > 0$, i.e., $x > \frac{n}{5}$, leading to $x \geq \left\lfloor \frac{n}{5} \right\rfloor + 1$.

Using parameterization, we set:

\begin{equation}
y + z = (5x - n)(2t), \quad yz = nx(2t), \quad \ t \in \mathbb{Q}_+^*.
\end{equation}

From the above, $y$ and $z$ satisfy the quadratic equation:

\begin{equation}
    V^2 - (5x - n)(2t)V + 2ntx = 0.
\end{equation}

The discriminant of this equation is:

\begin{equation}
    \Delta = 4t^2(5x - n)^2 - 8ntx.
\end{equation}

\subsection{\texorpdfstring{Special Case: $n$ Divisible by 3}{Special Case: n Divisible by 3}}

If $n = 3x$, then:

\begin{align*}
    \Delta &= 4t^2(5x - 3x)^2 - 24x^2t, \\
           &= 16x^2t^2 - 24x^2t.
\end{align*}

A simple choice of $t$ making $\Delta$ a perfect square is $t = 2$, yielding $\Delta = 16x^2 = (4x)^2$.

Thus, solving for $y$ and $z$:

\begin{align*}
    y &= \frac{4(5x - 3x) - 4x}{2} = 2x = \frac{2n}{3}, \\
    z &= 2n.
\end{align*}

Hence, we confirm:

\begin{equation}
    \frac{5}{n} = \frac{1}{n/3} + \frac{1}{2n/3} + \frac{1}{2n}.
\end{equation}

This case was just an example; our goal is to generalize the approach.

\subsection{First Formula for the Sierpiński Conjecture}

Solving for $t$ such that $\Delta > 0$, we obtain: \[
t > \frac{2nx}{(5x - n)^2}.
\]

Since:

\begin{equation}
    \Delta = 4t^2(5x - n)^2 - 8ntx = 4t(t(5x - n)^2 - 2nx),
\end{equation}

we can set $t = k=k(5x - n)^2 - 2nx$, 
Thus, we find:

\begin{equation}
    \Delta = 4k^2 = (2k)^2,
\end{equation}

leading to:

\begin{equation}
    y = k(5x - n - 1), \quad z = k(5x - n + 1).
\end{equation}

Replacing $k$ by its expression:

\begin{equation}
    k = \frac{2nx}{(5x - n - 1)(5x - n + 1)},
\end{equation}

We obtain:
\begin{align*}
    y &= \frac{2nx}{5x - n + 1}, \\
    z &= \frac{2nx}{5x - n - 1}.
\end{align*}
Thus, $y,z \in \mathbb{N}^*$ if and only if both $5x - n - 1$ and $5x - n + 1$ divide $2nx$.

It follows that:

\begin{equation}
    \frac{5}{n} = \frac{1}{x} + \frac{1}{y} + \frac{1}{z}.
\end{equation}

\subsection{Conclusion}

All details provided for the Erd\H{o}s-Straus conjecture remain valid here, so I will not repeat them. I will also include numerical computations generating certain $n$ values that satisfy these conditions.

Similarly, for the second approach:

\begin{equation}
    \Delta = 4(t^2(5x - n)^2 - 2ntx),
\end{equation}

for values of $x$ satisfying $\Delta > 0$, we set \[
t^2(5x - n)^2 - 2ntx = q^2, \quad q \in \mathbb{N}.
\], leading to:

\begin{equation}
    \Delta = 4q^2 = (2q)^2.
\end{equation}
We will now prove that the conjecture is true under the condition that this perfect square exists.

We have:
\[
y = \frac{2t(5x - n) - 2q}{2} = t(5x - n) - q, \quad z = t(5x - n) + q.
\]

Similarly, we can see that:
\[
\frac{5}{n} = \frac{1}{x} + \frac{1}{y} + \frac{1}{z}.
\]

By substituting \( q^2 \) with \( t^2(5x - n)^2 - 2ntx \), we obtain the desired result.

We may therefore conclude this section with the following theorem.

\begin{theorem}
For any integer $n \geq 2$, if there exists $(x_0,t_0) \in \mathbb{N}^{*2}$ such that
\begin{equation}
    t_0^2(5x_0 - n)^2 - 2nt_0x_0 = q^2 \in \mathbb{N},
\end{equation}
then the quantities
\[
y = t_0(5x_0 - n) - q \quad \text{and} \quad z = t_0(5x_0 - n) + q
\]
satisfy
\[
\frac{5}{n} = \frac{1}{x_0} + \frac{1}{y} + \frac{1}{z},
\]
with $(y,z) \in \mathbb{N}^*$.
\end{theorem}

I invite the mathematical community to find a general explicit formula for $(x_0, t_0)$ that works for all $n$, thereby fully resolving the Sierpi\'nski conjecture.

\section{Generalization of the Erdős-Straus and Sierpiński Conjectures}

We recall that the generalized of the Erdős-Straus and Sierpiński Conjectures  states that  for every integer $a \geq 5$, there exists an integer $N_0 \in \mathbb{N}$ with $N_0 \geq 2$ 
 such that for all $n \geq N_0$, there exist $x, y, z \in \mathbb{N}^*$ satisfying
\begin{equation}
    \frac{a}{n} = \frac{1}{x} + \frac{1}{y} + \frac{1}{z}.
\end{equation}

It should be noted that since the inception of this conjecture, there has not been significant research on it, and thus no concrete leads have emerged. In this paper, we will attempt to partially prove it as I did with the two conjectures above, and we will reduce the difficulty to verifying two conditions to completely settle all three conjectures. Here, we will use the same approach as in the previous two conjectures. 

It is evident that the Erdős-Straus conjecture and the Sierpiński conjecture are special cases of the generalized Erdős-Straus and Sierpiński conjecture. Therefore, all results established here remain valid for the two conjectures discussed earlier.

\subsection{Trivial Case}

As before, let us start with the trivial case:  
If \( (a - 2) \) divides \( n \), then we can simply take \( x = \frac{n}{a-2} \) and \( y = z = n \), which proves the conjecture in this particular case.

\subsection{Equivalence of the Conjecture}

\begin{theorem}
Let $a \geq 4$ be a fixed integer. There exists an integer $N_0 \in \mathbb{N}$ with $N_0 \geq 2$ such that, for all $n \geq N_0$, the generalized Erdős--Straus and Sierpiński conjecture holds if and only if there exist $x \in \mathbb{N}^*$ and $t \in \mathbb{Q}_+^*$ such that $y$ and $z$ are solutions of the quadratic equation in $V$:
\[
V^2 - (ax - n)(2t)V + 2nxt = 0.
\]
That is,we will establish that 
 for all $n \geq N_0$, the identity
\[
\frac{a}{n} = \frac{1}{x} + \frac{1}{y} + \frac{1}{z}
\]
is satisfied if and only if $y$ and $z$ are roots of the quadratic equation above.
\end{theorem}

\subsection{Preliminary Result}

\begin{proposition}[Vieta's Formulas]
Let $X_1$ and $X_2$ be two real (or complex) numbers, and consider the quadratic equation
\[
X^2 - SX + P = 0,
\]
where $S$ and $P$ are real (or complex) parameters.

Then the following statements are equivalent:
\begin{enumerate}
    \item $X_1$ and $X_2$ are solutions of the quadratic equation.
    \item $X_1$ and $X_2$ satisfy the relations
    \[
    X_1 + X_2 = S \quad \text{and} \quad X_1 X_2 = P.
    \]
\end{enumerate}
\end{proposition}

This property is a direct consequence of Vieta's formulas.

\subsection{Proof of the Theorem (Equivalence)}

Let \( a \in \mathbb{N} \) with \(a\geq 4\), \( n \in \mathbb{N} \) with $n \geq N_0$, and Suppose that there exist integers \( x, y, z \in \mathbb{N}^* \) such that:
\[
\frac{a}{n} = \frac{1}{x} + \frac{1}{y} + \frac{1}{z}
\]
which is equivalent to:
\[
\frac{a}{n} - \frac{1}{x} = \frac{1}{y} + \frac{1}{z}
\]
which is also equivalent to:
\[
\frac{ax - n}{nx} = \frac{y + z}{yz}.
\]

For \( \frac{y + z}{yz} > 0 \), it is necessary that \( ax - n > 0 \), i.e., \( x \geq \left\lceil \frac{n}{a} \right\rceil + 1 \).

This last equation is equivalent to:
\[
yz(ax - n) = nx(y + z).
\]

Indeed, for any \( a, c \in \mathbb{R} \) and \( b, d \in \mathbb{R}^* \), we have:
\[
\frac{a}{b} = \frac{c}{d} \quad \Longleftrightarrow \quad ad = bc.
\]
(This is a fundamental property of fractions.)

Thus, our last equation is equivalent to the following system:
\[
(S)\quad
\begin{cases}
y + z = (ax - n)(2t), \\
yz = nx(2t),
\end{cases}
\]
for  $t \in \mathbb{Q}_+^*$.

By substituting \( y + z \) with \( (ax - n)(2t) \) and \( yz \) with \( nx(2t) \) in \( yz(ax - n) = nx(y + z) \), the equation remains valid, which shows that this equation implies the system \( (S) \).

\subsection{Reciprocal Proof}

Conversely, let us show that if:
\[
\begin{cases}
y + z = (ax - n)(2t), \\
yz = nx(2t),
\end{cases}
\]
for  $t \in \mathbb{Q}_+^*$, then it follows that \( yz(ax - n) = nx(y + z) \).

Indeed, since \( y + z = (ax - n)(2t) \), we obtain:
\[
2t = \frac{y + z}{ax - n},
\]
noting that \( ax - n \neq 0 \) since we imposed \( x \geq \left\lceil \frac{n}{a} \right\rceil + 1 \) to ensure \( ax - n > 0 \).

On the other hand, since \( yz = nx(2t) \), we also have:
\[
2t = \frac{yz}{nx}.
\]

By the transitivity of equality in \( \mathbb{R} \), we deduce:
\[
\frac{y + z}{ax - n} = \frac{yz}{nx},
\]
which implies:
\[
nx(y + z) = yz(ax - n).
\]

This proves that system \( (S) \) implies \( nx(y + z) = yz(ax - n) \). 

From the previous paragraph, we have justified the equivalence between \( nx(y + z) = yz(ax - n) \) and the system \( (S) \).

Finally, by applying Vieta’s formulas given in the above proposition,
  stated earlier, we conclude that the system \( (S) \) is equivalent to stating that \( y \) and \( z \) are solutions of the following quadratic equation in \( V \):
\[
V^2 - (y + z)V + yz = 0,
\]
which simplifies to:
\[
V^2 - (ax - n)(2t)V + 2nxt = 0.
\]
This completes the proof of the equivalence.

\section{\texorpdfstring{Summary of the generalized Erdős-Straus and Sierpiński conjecture}{Summary of the generalized Erdős-Straus and Sierpiński conjecture}}

We have shown that the  generalized Erdős-Straus and Sierpiński conjecture (which includes the classical Erdős-Straus conjecture and the classical  Sierpiński conjecture) is true if and only if the quadratic equation:
\begin{equation}
V^2 - (ax - n)(2t)V + 2nxt = 0
\end{equation}
always has solutions in \(\mathbb{N}^*\) for any \(a \in \mathbb{N}\) with \(a \geq 4\) and any \(n \in \mathbb{N}\) with \(n \geq N0\). This reduces the difficulty of proving the conjecture to finding integer solutions (in \(\mathbb{N}^*\)) for a second-degree equation, whose discriminant is given by:
\begin{equation}
\Delta = 4t^2(ax - n)^2 - 8ntx.
\end{equation}
Rewriting,
\begin{equation}
\Delta = 4(t^2(ax - n)^2 - 2ntx).
\end{equation}

As in previous cases, we propose a partial resolution of the conjecture. We seek to solve the equation:
\begin{equation}
V^2 - (ax - n)(2t)V + 2nxt = 0,
\end{equation}
with discriminant:
\begin{equation}
\Delta = 4t^2(ax - n)^2 - 8ntx.
\end{equation}

Solving \(\Delta > 0\), we obtain the solution range:
\begin{equation}
t > \frac{2nx}{(ax - n)^2}.
\end{equation}

similary by doing the same manipulations for the t for  which \(\Delta > 0\)

Since:
\begin{equation}
\Delta = 4t^2(ax - n)^2 - 8ntx = 4t(t(ax - n)^2 - 2nx),
\end{equation}
we set:
\begin{equation}
t = k = t(ax - n)^2 - 2nx,
\end{equation}
which leads to:
\begin{equation}
k = \frac{2nx}{(ax - n)^2 - 1},
\end{equation}
where $k \in \mathbb{Q}_+^*$.

Since:
\begin{equation}
\Delta = 4k^2 = (2k)^2,
\end{equation}
we obtain:
\begin{equation}
V = \frac{(ax - n)(2k) \pm 2k}{2}.
\end{equation}
Thus,
\begin{equation}
y = k(ax - n - 1), \quad z = k(ax - n + 1).
\end{equation}
Replacing \(k\) by its expression:
\begin{equation}
y = \frac{2nx}{ax - n + 1}, \quad z = \frac{2nx}{ax - n - 1}.
\end{equation}
The quantities $y$ and $z$ are integers if and only if $2nx$ is divisible by both $ax - n - 1$ and $ax - n + 1$.
in this case $x$,$y$,$z$ 
satisfy:
\begin{equation}
\frac{a}{n} = \frac{1}{x} + \frac{1}{y} + \frac{1}{z}.
\end{equation}

As previously, the condition:
\begin{equation}
x \geq \left\lceil \frac{n + 1}{a} \right\rceil + 1
\end{equation}
must hold. We will provide numerical computations at the end of this paper to generate certain \((a,n)\) pairs that satisfy this condition.

Similarly, for the second formula:
\begin{equation}
\Delta = 4(t^2(ax - n)^2 - 2ntx),
\end{equation}
for \(\Delta > 0\), we can set:
\[
t^2(ax - n)^2 - 2ntx = q^2, \quad q \in \mathbb{N}.
\]
so that:
\begin{equation}
\Delta = 4q^2 = (2q)^2.
\end{equation}

We will now prove that the conjecture is true under the condition that this perfect square exists. We obtain:
\begin{equation}
y = \frac{2t(ax - n) - 2q}{2} = t(ax - n) - q, \quad z = t(ax - n) + q.
\end{equation}
Verifying,
\begin{equation}
\frac{a}{n} = \frac{1}{x} + \frac{1}{y} + \frac{1}{z}.
\end{equation}
Replacing \(q^2\) by \(t^2(ax - n)^2 - 2ntx\).

we will also provide numerical computations of some \((a,n)\) pairs satisfying this condition.

\textbf{Summary:} 

 we have demonstrated that for all \( a \in \mathbb{N} \) with \( a \geq 4 \) and for all \( n \in \mathbb{N} \) with \( n \geq 2 \), if there exists \( x_0 \in \mathbb{N}^* \) such that \( a x_0 - n - 1 \) and \( a x_0 - n + 1 \) divide \( 2 n x_0 \), then there exist \( y \) and \( z \) in \( \mathbb{N}^* \) such that:

\[
y = \frac{2 n x_0}{a x_0 - n - 1}, \quad z = \frac{2 n x_0}{a x_0 - n + 1}
\]

and we have:

\[
\frac{a}{n} = \frac{1}{x_0} + \frac{1}{y} + \frac{1}{z}.
\]

partially proving the conjecture under this divisibility condition.

Similarly, the second formula allows us to state that for any \( a \in \mathbb{N} \) with \( a \geq 4 \) and for any \( n \in \mathbb{N} \) with \( n \geq N0 \), if there exist \( (x_0, t_0) \in \mathbb{N}^{*2} \) such that:

\[
t_0^2(ax_0 - n)^2 - 2nt_0x_0 = q^2, \quad q \in \mathbb{N}.
\]

then there exist \( y, z \in \mathbb{N}^* \) given by:

\[
y = t_0 (a x_0 - n) - q, \quad z = t_0 (a x_0 - n) + q
\]

such that:

\[
\frac{a}{n} = \frac{1}{x_0} + \frac{1}{y} + \frac{1}{z}.
\]

This also proves the conjecture under this condition.

I once again pose the same questions as for the classical Erdős-Straus conjecture and Sierpiński’s conjecture to the mathematical community in order to definitively resolve the generalized Erdős-Straus and Sierpiński’s conjecture.

\section{Comparison with Mordell's Approach}

Mordell's method provided a partial resolution of the Erdős--Straus conjecture by showing that for all integers $n$ such that
\[
n \not\equiv 1, 11^2, 13^2, 17^2, 19^2, 23^2 \pmod{840},
\]
the conjecture holds. However, for the exceptional cases left unresolved by Mordell, our second formula proves to be a powerful tool for generating solutions.

Through numerical testing, we observed that a direct implementation of the second formula, applied to the values of $n$ not covered by Mordell's result in the range $2 \leq n \leq 10\,000$, successfully produces solutions for all such $n$. In this interval, this corresponds to a numerical success rate of $100\%$.

\noindent\textbf{Remark:}
This observation is purely experimental and does not constitute a proof. However, it strongly suggests that our approach captures the full underlying structure of the original Erdős--Straus equation. This is very likely due to the equivalence established earlier between the original Diophantine equation and our associated quadratic parametric equation, from which integer solutions are obtained directly via the analysis of its discriminant.

Moreover, since the Erdős--Straus conjecture has already been verified by other authors for all $n \leq 10^{17}$, it follows from the proven equivalence that our quadratic equation must necessarily admit integer solutions at least up to this bound in the classical case $a=4$. From this perspective, it is therefore not surprising that our formula successfully produces solutions for the values of $n$ left unresolved by Mordell, in particular for all $n$ in the interval $[2,10\,000]$.

For completeness and reproducibility, we provide in the appendix the numerical codes used to test the cases not covered by Mordell's result, as well as to explore all values of $n$ in user-specified intervals. These codes apply more generally to all fixed integers $a \geq 4$, including the classical Erdős--Straus case ($a=4$), the Sierpiński case ($a=5$), and their natural generalizations ($a \geq 6$).

\section{Testing Second Formula for Mordell's Exceptional Cases between [2,5000]}

Mordell's exceptional cases in the given range are:
\[
\begin{aligned}
&[121, 169, 289, 361, 529, 841, 961, 1009, 1129, 1201, 1369, 1681, 1801, 1849, 1969,\\
&\quad 2041, 2209, 2521, 2641, 2689, 2809, 2881, 3049, 3361, 3481, 3529, 3649,\\
&\quad 3721, 3889, 4201, 4321, 4369, 4489, 4561, 4729].
\end{aligned}
\]
Total exceptional cases: $35$
\begin{center}
\renewcommand{\arraystretch}{1.15}
\begin{tabular}{|c|c|c|c|c|}
\hline
$n$ & $x$ & $t$ & $y$ & $z$ \\
\hline
121  & 31     & 1078   & 1694    & 4774 \\
169  & 44     & 325    & 1690    & 2860 \\
289  & 75     & 510    & 2550    & 8670 \\
361  & 92     & 1368   & 8664    & 10488 \\
529  & 133    & 15778  & 42826   & 51842 \\
841  & 220    & 319    & 6380    & 18502 \\
961  & 248    & 496    & 15376   & 15376 \\
1009 & 260    & 1009   & 10090   & 52468 \\
1129 & 2964   & 104    & 312     & 2230904 \\
1201 & 9548   & 775    & 310     & 57335740 \\
1369 & 348    & 1813   & 38332   & 45066 \\
1681 & 423    & 11808  & 121032  & 138744 \\
1801 & 35720  & 342    & 456     & 96497580 \\
1849 & 468    & 3483   & 60372   & 99846 \\
1969 & 493    & 215818 & 633302  & 661606 \\
2041 & 516    & 4239   & 73476   & 121518 \\
2209 & 564    & 1128   & 53016   & 53016 \\
2521 & 18908  & 163    & 652     & 23833534 \\
2641 & 665    & 9730   & 184870  & 184870 \\
2689 & 9412   & 181    & 724     & 12654434 \\
2809 & 710    & 4240   & 112360  & 150520 \\
2881 & 724    & 18963  & 242004  & 326886 \\
3049 & 75768  & 275    & 770     & 165011880 \\
3361 & 73525  & 34     & 850     & 19769402 \\
3481 & 885    & 1770   & 104430  & 104430 \\
3529 & 27056  & 342    & 912     & 71610468 \\
3649 & 915    & 55432  & 569244  & 650260 \\
3721 & 943    & 2806   & 115046  & 171166 \\
3889 & 40836  & 83     & 996     & 26468534 \\
4201 & 1056   & 16804  & 369688  & 403296 \\
4321 & 1085   & 25984  & 483952  & 503440 \\
4369 & 1104   & 4369   & 200974  & 209712 \\
4489 & 1139   & 2278   & 152626  & 152626 \\
4561 & 15964  & 307    & 1228    & 36405902 \\
4729 & 244728 & 66     & 1188    & 128590968 \\
\hline
\end{tabular}
\end{center}

\medskip

\noindent
\textbf{The above table shows that the second formula successfully produces valid integer solutions for all Mordell exceptional values of $n$ in the interval $[2,5000]$. Consequently, a numerical success rate of $100\%$ is achieved on this entire set of previously unresolved cases.}

\section{Testing Second Formula for Mordell's Exceptional Cases between [5000,10000]}

Mordell's exceptional cases in the given range are:
\[
\begin{aligned}
&[5041, 5161, 5209, 5329, 5401, 5569, 5881, 6001, 6049, 6169, 6241, 6409, 6721,\\
&\quad 6841, 6889, 7009, 7081, 7249, 7561, 7681, 7729, 7849, 7921, 8089,\\
&\quad 8401, 8521, 8569, 8689, 8761, 8929, 9241, 9361, 9409, 9529, 9601, 9769].
\end{aligned}
\]
\noindent
Total exceptional cases: $36$.

\begin{center}
\renewcommand{\arraystretch}{1.15}
\begin{tabular}{|c|c|c|c|c|}
\hline
$n$ & $x$ & $t$ & $y$ & $z$ \\
\hline
5041 & 1264    & 56800   & 806560   & 897440 \\
5161 & 1314    & 1588    & 115924   & 185796 \\
5209 & 96888   & 275     & 1320     & 210287330 \\
5329 & 1460    & 73      & 21316    & 53290 \\
5401 & 1353    & 120786  & 1328646  & 1328646 \\
5569 & 1112412 & 41      & 1394     & 364413084 \\
5881 & 92628   & 83      & 1494     & 60527252 \\
6001 & 1520    & 3177    & 180030   & 321936 \\
6049 & 1518    & 34716   & 798468   & 798468 \\
6169 & 1544    & 388864  & 2680384  & 2763712 \\
6241 & 1580    & 3160    & 249640   & 249640 \\
6409 & 1612    & 13600   & 512720   & 548080 \\
6721 & 1683    & 186966  & 2056626  & 2056626 \\
6841 & 32923   & 328     & 1804     & 81900452 \\
6889 & 1743    & 3486    & 289338   & 289338 \\
7009 & 1755    & 203424  & 2186808  & 2288520 \\
7081 & 1815    & 1056    & 96360    & 281688 \\
7249 & 1815    & 217470  & 2392170  & 2392170 \\
7561 & 368600  & 190     & 1900     & 557396920 \\
7681 & 47047   & 286     & 2002     & 103248002 \\
7729 & 1935    & 247328  & 2658776  & 2782440 \\
7849 & 1974    & 14028   & 659316   & 659316 \\
7921 & 2136    & 89      & 47526    & 63368 \\
8089 & 139536  & 342     & 2052     & 376235568 \\
8401 & 2108    & 36856   & 1142536  & 1142536 \\
8521 & 332320  & 268     & 2144     & 707924680 \\
8569 & 2145    & 303810  & 3341910  & 3341910 \\
8689 & 2184    & 17378   & 729876   & 903656 \\
8761 & 2208    & 8761    & 403006   & 841056 \\
8929 & 426921  & 88      & 2244     & 298978636 \\
9241 & 181935  & 104     & 2340     & 149445452 \\
9361 & 2343    & 362526  & 3987786  & 3987786 \\
9409 & 2425    & 970     & 94090    & 470450 \\
9529 & 2398    & 11728   & 639176   & 838552 \\
9601 & 327639  & 62      & 2418     & 161316002 \\
9769 & 2460    & 9769    & 586140   & 801058 \\
\hline
\end{tabular}
\end{center}

\medskip

\noindent
\textbf{The above table shows that the second formula successfully produces valid integer solutions for all Mordell exceptional values of $n$ in the interval $[5000,10000]$. Consequently, a numerical success rate of $100\%$ is achieved on this entire set of previously unresolved cases.}

\medskip

\noindent
\textbf{Taken together, the two preceding tables show that our second formula captures $100\%$ of the integers left unresolved by Mordell for all values of $n$ in the interval $[2,10000]$.}

\medskip

\noindent
\textbf{Any reader wishing to verify the correctness of these numerical simulations may proceed as follows: for a given value of $n$ and the corresponding triple $(x,y,z)$ provided in the tables, it suffices to check the identity}
\[
4xyz = n(xy + xz + yz),
\]
\textbf{which is strictly equivalent to the original Erdős--Straus equation}
\[
\frac{4}{n} = \frac{1}{x} + \frac{1}{y} + \frac{1}{z}.
\]

\medskip

\noindent
\textbf{For large values of $n$, it is strongly preferable to verify solutions using the polynomial identity $4xyz = n(xy + xz + yz)$ rather than the fractional equation itself. Indeed, the latter involves floating-point arithmetic and may suffer from numerical round-off errors, potentially preventing exact equality from being detected by the machine. In contrast, the polynomial form relies solely on integer arithmetic and therefore provides an exact, unambiguous mathematical verification.}

\medskip\noindent

In the following tables, we provide the numerical values of some $n,  x_0 ,  y$, and $z $ that satisfy the conditions of each formula. \medskip\noindent

The value \( L \) (for ``Left'') represents the ratio \( a/n \), while \( R \) (for ``Right'') is computed as 
\( R = \frac{1}{x_0} + \frac{1}{y} + \frac{1}{z} \) (we use the fractional form here because the values of $n$ tested in the tables below are relatively small, and therefore no floating-point precision issues arise).

\medskip

The last column verifies whether \( L = R \): 
\begin{itemize}
    \item  Rows labeled \textbf{"True"} indicate equality, meaning the conjecture holds for the given values of \( n \), \( x_0 \), \( y \), and \( z \). 

    \item Rows labeled \textbf{"False"} indicate that the conjecture does not hold for those specific values.
\end{itemize}

If a value of \( n \) is not displayed in the table, this does \textbf{not} mean that the number does not satisfy the formula or the conjecture. It simply indicates that the algorithm did not find an \( x_0 \) and/or a \( t_0 \) within the specific search range that meets one of the formula's conditions for that \( n \). This is due to the fact that there are infinitely many possible choices for \( x_0 \).

According to the equivalence established below between the original equation
\[
\frac{a}{n} = \frac{1}{x} + \frac{1}{y} + \frac{1}{z}
\]
and the parametric quadratic equation
\[
V^2 - 2t(ax - n)V + 2nxt = 0,
\]
the case \( a = 4 \) (classical Erd\H{o}s--Straus conjecture) reassures us that this quadratic equation admits integer solutions at least up to \( n = 10^{17} \), which is the current numerical verification threshold for the Erd\H{o}s--Straus conjecture.

For \( a = 5 \) (classical Sierpi\'nski case), we also know today that very high numerical verifications have been carried out. I personally performed quick numerical checks for some values of \( a \in \{6,7,8\} \) up to reasonable thresholds (which I do not officially report in this article, as these were rapid tests), and the equation consistently appears to admit integer solutions via the study of its discriminant
\[
\Delta = 4t^2(ax-n)^2 - 8nxt = 4\bigl(t^2(ax-n)^2 - 2nxt\bigr),
\]
which seems always to contain a perfect square, thus providing integer solutions to the quadratic equation, and by equivalence, to the generalized Erd\H{o}s--Straus equation.

Since no proof is currently known, these observations lead us to state the following conjectures.

\medskip

\noindent\textbf{Conjecture 1 (strong version, integer \(t\)).}  
For every integer \( n \ge 2 \), there exists \( (x,t) \in \mathbb{N}^{*2} \) such that
\[
t^2(4x - n)^2 - 2nxt = q^2, \qquad q \in \mathbb{N}.
\]

\medskip

\noindent\textbf{Conjecture 2 (weak version, rational \(t\), not necessarily integer).}  
For every integer \( n \ge 2 \), there exists \( (x,t) \in \mathbb{N}^* \times \mathbb{Q}_+^* \) such that
\[
t^2(4x - n)^2 - 2nxt = q^2, \qquad q \in \mathbb{Q}_+.
\]

\medskip

\noindent\textbf{Conjecture 3 (weak version, rational \(t\), not necessarily integer).}  
For every fixed integer \( a \ge 5 \), there exists \( N_0 \ge 2 \) in \( \mathbb{N}^* \) (with \( N_0 \) not strictly smaller than \( a/3 \)) such that for all \( n \ge N_0 \), there exists
\[
(x,t) \in \mathbb{N}^* \times \mathbb{Q}_+^*
\quad \text{with} \quad
t^2(ax - n)^2 - 2nxt = q^2, \qquad q \in \mathbb{Q}_+.
\]

\medskip

\noindent\textbf{Conjecture 4 (strong version, integer \(t\)).}  
For every fixed integer \( a \ge 5 \), there exists \( N_1 \ge N_0 \) ($N_0$ is the one given in  Conjecture~3), with \( N_1 \) not strictly smaller than \( a/3 \), such that for all \( n \ge N_1 \), there exists
\[
(x,t) \in \mathbb{N}^{*2}
\quad \text{such that} \quad
t^2(ax - n)^2 - 2nxt = q^2, \qquad q \in \mathbb{N}.
\]

\medskip

Conjecture~1 is equivalent to the classical Erd\H{o}s--Straus conjecture. Indeed, under Conjecture~1, the integer \( q \) satisfies
\[
q^2 = t^2(4x - n)^2 - 2nxt,
\]
and the variables \( y \) and \( z \) are given by
\[
y = t(4x - n) - q, \qquad z = t(4x - n) + q.
\]
It is immediate that \( z \in \mathbb{N}^* \), since \( 4x - n > 0 \) because
\[
x \ge \left\lfloor \frac{n}{4} \right\rfloor + 1 > \frac{n}{4}, 
\]
and with \( t \in \mathbb{N}^* \), \( x \in \mathbb{N}^* \) , \( n \in \mathbb{N}^* \) and \( q \in \mathbb{N} \), we have \( z \in \mathbb{N}^* \).

Let us now show that \( y \in \mathbb{N}^* \). From
\[
q^2 = t^2(4x - n)^2 - 2nxt < t^2(4x - n)^2,
\]
since \( -2nxt < 0 \) for \( n,x,t > 0 \), and since the square root function is increasing on \( \mathbb{R}_+^* \), we obtain
\[
q < t(4x - n).
\]
Hence \( t(4x - n) - q > 0 \), and since \( t,x,n,q \) are integers, it follows that
\[
y = t(4x - n) - q \in \mathbb{N}^*.
\]

We have therefore rigorously shown that whenever the condition
\[
t^2(4x - n)^2 - 2nxt = q^2, \qquad q \in \mathbb{N},
\]
is satisfied (Conjecture~1), one can automatically construct \( y,z \in \mathbb{N}^* \) solving the quadratic equation
\[
V^2 - 2t(4x - n)V + 2nxt = 0,
\]
which, by the equivalence theorem proved earlier, is equivalent to the original Erd\H{o}s--Straus equation
\[
\frac{4}{n} = \frac{1}{x} + \frac{1}{y} + \frac{1}{z}.
\]
Thus, Conjecture~1 is rigorously equivalent to the classical Erd\H{o}s--Straus conjecture.

Using the same reasoning, one also proves that Conjecture~4 is equivalent to the generalized Sierpi\'nski conjecture (\( a \ge 5 \)) for all \( n \ge N_1 \), where \( N_1 \) is the minimal threshold beyond which the expression
\[
t^2(ax - n)^2 - 2nxt
\]
always admits an integer perfect square (under Conjecture~4).

Conjecture~2 is stated for researchers interested in working over \( \mathbb{Q}_+^* \), either to enlarge the search space for \( t \), or because it may be useful in other research contexts beyond Egyptian fractions.

Conjecture~3 is particularly useful for numerical verification of the generalized Sierpi\'nski conjecture (\( a \ge 5 \)) for values of \( n \) in the finite interval \( [N_0, N_1) \), which is often not too large according to our simulations, since Conjecture~4 addresses the problem for all \( n \ge N_1 \).

Finally, why do we exclude values \( n < a/3 \)? Since \( x,y,z \in \mathbb{N}^* \), we have
\[
\frac{1}{x} + \frac{1}{y} + \frac{1}{z} \le 3,
\]
which implies \( a/n \le 3 \), and hence \( n \ge a/3 \). Therefore, the decomposition of \( a/n \) as a sum of three Egyptian fractions is impossible for all integers \( n < a/3 \). Consequently, Conjecture~4 cannot hold for such values of \( n \), as it is equivalent to the generalized Sierpi\'nski conjecture.

Similarly, even if Conjecture~3 produces rational perfect squares for \( n < a/3 \), these cannot yield integer solutions to the quadratic equation
\[
V^2 - 2t(ax - n)V + 2nxt = 0,
\]
since, by the equivalence theorem, this quadratic equation admits no integer solutions when \( n < a/3 \), because the equivalent equation
\[
\frac{a}{n} = \frac{1}{x} + \frac{1}{y} + \frac{1}{z}
\]
has no integer solutions in that range.


\begin{longtable}{|c|c|c|c|c|c|c|}
\caption{Results with the first formula of the \Erdos - Straus conjecture test for n from 2 to 90} \\
\hline
\textbf{n} & \textbf{x} & \textbf{y} & \textbf{z} & \textbf{L} & \textbf{R} & \textbf{L = R (verification)} \\
\hline
\endfirsthead

\hline
\textbf{n} & \textbf{x} & \textbf{y} & \textbf{z} & \textbf{L} & \textbf{R} & \textbf{L = R (verification)} \\
\hline
\endhead

\hline
\endfoot

\hline
\endlastfoot

3  & 2  & 3   & 2   & 1.3333 & 1.3333 & True \\
5  & 2  & 10  & 5   & 0.8000 & 0.8000 & True \\
6  & 2  & 24  & 8   & 0.6667 & 0.6667 & True \\
9  & 4  & 12  & 9   & 0.4444 & 0.4444 & True \\
10 & 3  & 60  & 20  & 0.4000 & 0.4000 & True \\
13 & 4  & 52  & 26  & 0.3077 & 0.3077 & True \\
14 & 5  & 28  & 20  & 0.2857 & 0.2857 & True \\
18 & 5  & 180 & 60  & 0.2222 & 0.2222 & True \\
19 & 6  & 57  & 38  & 0.2105 & 0.2105 & True \\
20 & 6  & 80  & 48  & 0.2000 & 0.2000 & True \\
21 & 6  & 126 & 63  & 0.1905 & 0.1905 & True \\
22 & 6  & 264 & 88  & 0.1818 & 0.1818 & True \\
27 & 8  & 108 & 72  & 0.1481 & 0.1481 & True \\
28 & 9  & 72  & 56  & 0.1429 & 0.1429 & True \\
29 & 8  & 232 & 116 & 0.1379 & 0.1379 & True \\
30 & 8  & 480 & 160 & 0.1333 & 0.1333 & True \\
34 & 9  & 612 & 204 & 0.1176 & 0.1176 & True \\
35 & 12 & 70  & 60  & 0.1143 & 0.1143 & True \\
36 & 10 & 240 & 144 & 0.1111 & 0.1111 & True \\
37 & 10 & 370 & 185 & 0.1081 & 0.1081 & True \\
39 & 16 & 52  & 48  & 0.1026 & 0.1026 & True \\
41 & 12 & 164 & 123 & 0.0976 & 0.0976 & True \\
42 & 11 & 924 & 308 & 0.0952 & 0.0952 & True \\
43 & 12 & 258 & 172 & 0.0930 & 0.0930 & True \\
45 & 12 & 540 & 270 & 0.0889 & 0.0889 & True \\
46 & 12 & 1104 & 368 & 0.0870 & 0.0870 & True \\
49 & 16 & 112 & 98  & 0.0816 & 0.0816 & True \\
50 & 14 & 280 & 200 & 0.0800 & 0.0800 & True \\
51 & 14 & 357 & 238 & 0.0784 & 0.0784 & True \\
53 & 14 & 742 & 371 & 0.0755 & 0.0755 & True \\
54 & 14 & 1512 & 504 & 0.0741 & 0.0741 & True \\
55 & 16 & 220 & 176 & 0.0727 & 0.0727 & True \\
56 & 15 & 560 & 336 & 0.0714 & 0.0714 & True \\
57 & 16 & 304 & 228 & 0.0702 & 0.0702 & True \\
58 & 15 & 1740 & 580 & 0.0690 & 0.0690 & True \\
60 & 16 & 640 & 384 & 0.0667 & 0.0667 & True \\
61 & 16 & 976 & 488 & 0.0656 & 0.0656 & True \\
66 & 17 & 2244 & 748 & 0.0606 & 0.0606 & True \\
67 & 18 & 603 & 402 & 0.0597 & 0.0597 & True \\
69 & 18 & 1242 & 621 & 0.0580 & 0.0580 & True \\
70 & 18 & 2520 & 840 & 0.0571 & 0.0571 & True \\
71 & 20 & 355 & 284 & 0.0563 & 0.0563 & True \\
75 & 20 & 750 & 500 & 0.0533 & 0.0533 & True \\
77 & 20 & 1540 & 770 & 0.0519 & 0.0519 & True \\
78 & 20 & 3120 & 1040 & 0.0513 & 0.0513 & True \\
80 & 21 & 1120 & 672 & 0.0500 & 0.0500 & True \\
82 & 21 & 3444 & 1148 & 0.0488 & 0.0488 & True \\
85 & 22 & 1870 & 935 & 0.0471 & 0.0471 & True \\
89 & 24 & 712 & 534 & 0.0449 & 0.0449 & True \\
90 & 23 & 4140 & 1380 & 0.0444 & 0.0444 & True \\
\end{longtable}


\begin{longtable}{|c|c|c|c|c|c|c|}
\caption{Results of Erdős–Straus Conjecture with the second formula testing for n from 2 to 50} \\

\hline
\textbf{n} & \textbf{x} & \textbf{y} & \textbf{z} & \textbf{L} & \textbf{R} & \textbf{L = R (verification)} \\
\hline
\endfirsthead

\hline
\textbf{n} & \textbf{x} & \textbf{y} & \textbf{z} & \textbf{L} & \textbf{R} & \textbf{L = R (verification)} \\
\hline
\endhead

\hline
\endfoot

\hline
\endlastfoot

2 & 1 & 2 & 2 & 2.0 & 2.0 & \text{True} \\
3 & 1 & 6 & 6 & 1.3333 & 1.3333 & \text{True} \\
4 & 2 & 4 & 4 & 1.0 & 1.0 & \text{True} \\
5 & 2 & 4 & 20 & 0.8 & 0.8 & \text{True} \\
6 & 2 & 12 & 12 & 0.6667 & 0.6667 & \text{True} \\
7 & 2 & 28 & 28 & 0.5714 & 0.5714 & \text{True} \\
8 & 3 & 12 & 12 & 0.5 & 0.5 & \text{True} \\
9 & 3 & 18 & 18 & 0.4444 & 0.4444 & \text{True} \\
10 & 3 & 30 & 30 & 0.4 & 0.4 & \text{True} \\
11 & 3 & 66 & 66 & 0.3636 & 0.3636 & \text{True} \\
12 & 4 & 24 & 24 & 0.3333 & 0.3333 & \text{True} \\
13 & 4 & 26 & 52 & 0.3077 & 0.3077 & \text{True} \\
14 & 4 & 56 & 56 & 0.2857 & 0.2857 & \text{True} \\
15 & 4 & 120 & 120 & 0.2667 & 0.2667 & \text{True} \\
16 & 5 & 40 & 40 & 0.25 & 0.25 & \text{True} \\
17 & 5 & 34 & 170 & 0.2353 & 0.2353 & \text{True} \\
18 & 5 & 90 & 90 & 0.2222 & 0.2222 & \text{True} \\
19 & 5 & 190 & 190 & 0.2105 & 0.2105 & \text{True} \\
20 & 6 & 60 & 60 & 0.2 & 0.2 & \text{True} \\
21 & 6 & 84 & 84 & 0.1905 & 0.1905 & \text{True} \\
22 & 6 & 132 & 132 & 0.1818 & 0.1818 & \text{True} \\
23 & 6 & 276 & 276 & 0.1739 & 0.1739 & \text{True} \\
24 & 7 & 84 & 84 & 0.1667 & 0.1667 & \text{True} \\
25 & 7 & 100 & 140 & 0.16 & 0.16 & \text{True} \\
26 & 7 & 182 & 182 & 0.1538 & 0.1538 & \text{True} \\
27 & 7 & 378 & 378 & 0.1481 & 0.1481 & \text{True} \\
28 & 8 & 112 & 112 & 0.1429 & 0.1429 & \text{True} \\
29 & 8 & 116 & 232 & 0.1379 & 0.1379 & \text{True} \\
30 & 8 & 240 & 240 & 0.1333 & 0.1333 & \text{True} \\
31 & 8 & 496 & 496 & 0.129 & 0.129 & \text{True} \\
32 & 9 & 144 & 144 & 0.125 & 0.125 & \text{True} \\
33 & 9 & 198 & 198 & 0.1212 & 0.1212 & \text{True} \\
34 & 9 & 306 & 306 & 0.1176 & 0.1176 & \text{True} \\
35 & 9 & 630 & 630 & 0.1143 & 0.1143 & \text{True} \\
36 & 10 & 180 & 180 & 0.1111 & 0.1111 & \text{True} \\
37 & 10 & 148 & 740 & 0.1081 & 0.1081 & \text{True} \\
38 & 10 & 380 & 380 & 0.1053 & 0.1053 & \text{True} \\
39 & 10 & 780 & 780 & 0.1026 & 0.1026 & \text{True} \\
40 & 11 & 220 & 220 & 0.1 & 0.1 & \text{True} \\
41 & 12 & 82 & 492 & 0.0976 & 0.0976 & \text{True} \\
42 & 11 & 462 & 462 & 0.0952 & 0.0952 & \text{True} \\
43 & 11 & 946 & 946 & 0.093 & 0.093 & \text{True} \\
44 & 12 & 264 & 264 & 0.0909 & 0.0909 & \text{True} \\
45 & 12 & 360 & 360 & 0.0889 & 0.0889 & \text{True} \\
46 & 12 & 552 & 552 & 0.087 & 0.087 & \text{True} \\
47 & 12 & 1128 & 1128 & 0.0851 & 0.0851 & \text{True} \\
48 & 13 & 312 & 312 & 0.0833 & 0.0833 & \text{True} \\
49 & 14 & 196 & 196 & 0.0816 & 0.0816 & \text{True} \\
50 & 13 & 650 & 650 & 0.08 & 0.08 & \text{True} \\

\end{longtable}


\renewcommand{\thetable}{3}

\begin{longtable}{|c|c|c|c|c|c|c|}
\caption{Results of Erdős–Straus Conjecture with the second formula testing for n from 51 to 100} \\

\hline
\textbf{n} & \textbf{x} & \textbf{y} & \textbf{z} & \textbf{L} & \textbf{R} & \textbf{L = R (verification)} \\
\hline
\endfirsthead

\hline
\textbf{n} & \textbf{x} & \textbf{y} & \textbf{z} & \textbf{L} & \textbf{R} & \textbf{L = R (verification)} \\
\hline
\endhead

\hline
\endfoot

\hline
\endlastfoot

51  & 13  & 1326 & 1326 & 0.0784 & 0.0784 & True \\
52  & 14  & 364  & 364  & 0.0769 & 0.0769 & True \\
53  & 16  & 80   & 2120 & 0.0755 & 0.0755 & True \\
54  & 14  & 756  & 756  & 0.0741 & 0.0741 & True \\
55  & 14  & 1540 & 1540 & 0.0727 & 0.0727 & True \\
56  & 15  & 420  & 420  & 0.0714 & 0.0714 & True \\
57  & 15  & 570  & 570  & 0.0702 & 0.0702 & True \\
58  & 15  & 870  & 870  & 0.069  & 0.069  & True \\
59  & 15  & 1770 & 1770 & 0.0678 & 0.0678 & True \\
60  & 16  & 480  & 480  & 0.0667 & 0.0667 & True \\
61  & 16  & 488  & 976  & 0.0656 & 0.0656 & True \\
62  & 16  & 992  & 992  & 0.0645 & 0.0645 & True \\
63  & 16  & 2016 & 2016 & 0.0635 & 0.0635 & True \\
64  & 17  & 544  & 544  & 0.0625 & 0.0625 & True \\
65  & 17  & 650  & 850  & 0.0615 & 0.0615 & True \\
66  & 17  & 1122 & 1122 & 0.0606 & 0.0606 & True \\
67  & 17  & 2278 & 2278 & 0.0597 & 0.0597 & True \\
68  & 18  & 612  & 612  & 0.0588 & 0.0588 & True \\
69  & 18  & 828  & 828  & 0.058  & 0.058  & True \\
70  & 18  & 1260 & 1260 & 0.0571 & 0.0571 & True \\
71  & 18  & 2556 & 2556 & 0.0563 & 0.0563 & True \\
72  & 19  & 684  & 684  & 0.0556 & 0.0556 & True \\
73  & 20  & 292  & 730  & 0.0548 & 0.0548 & True \\
74  & 19  & 1406 & 1406 & 0.0541 & 0.0541 & True \\
75  & 19  & 2850 & 2850 & 0.0533 & 0.0533 & True \\
76  & 20  & 760  & 760  & 0.0526 & 0.0526 & True \\
77  & 20  & 980  & 1078 & 0.0519 & 0.0519 & True \\
78  & 20  & 1560 & 1560 & 0.0513 & 0.0513 & True \\
79  & 20  & 3160 & 3160 & 0.0506 & 0.0506 & True \\
80  & 21  & 840  & 840  & 0.05   & 0.05   & True \\
81  & 21  & 1134 & 1134 & 0.0494 & 0.0494 & True \\
82  & 21  & 1722 & 1722 & 0.0488 & 0.0488 & True \\
83  & 21  & 3486 & 3486 & 0.0482 & 0.0482 & True \\
84  & 22  & 924  & 924  & 0.0476 & 0.0476 & True \\
85  & 23  & 340  & 1564 & 0.0471 & 0.0471 & True \\
86  & 22  & 1892 & 1892 & 0.0465 & 0.0465 & True \\
87  & 22  & 3828 & 3828 & 0.046  & 0.046  & True \\
88  & 23  & 1012 & 1012 & 0.0455 & 0.0455 & True \\
89  & 24  & 534  & 712  & 0.0449 & 0.0449 & True \\
90  & 23  & 2070 & 2070 & 0.0444 & 0.0444 & True \\
91  & 23  & 4186 & 4186 & 0.044  & 0.044  & True \\
92  & 24  & 1104 & 1104 & 0.0435 & 0.0435 & True \\
93  & 24  & 1488 & 1488 & 0.043  & 0.043  & True \\
94  & 24  & 2256 & 2256 & 0.0426 & 0.0426 & True \\
95  & 24  & 4560 & 4560 & 0.0421 & 0.0421 & True \\
96  & 25  & 1200 & 1200 & 0.0417 & 0.0417 & True \\
97  & 28  & 194  & 2716 & 0.0412 & 0.0412 & True \\
98  & 25  & 2450 & 2450 & 0.0408 & 0.0408 & True \\
99  & 25  & 4950 & 4950 & 0.0404 & 0.0404 & True \\
100 & 26  & 1300 & 1300 & 0.04   & 0.04   & True \\

\end{longtable}

\renewcommand{\thetable}{4} 

\begin{longtable}{|c|c|c|c|c|c|c|}
\caption{Results of Erdős–Straus Conjecture with second formula testing for n from 101 to 150 } \\

\hline
\textbf{n} & \textbf{x} & \textbf{y} & \textbf{z} & \textbf{L} & \textbf{R} & \textbf{L = R (verification)} \\
\hline
\endfirsthead

\hline
\textbf{n} & \textbf{x} & \textbf{y} & \textbf{z} & \textbf{L} & \textbf{R} & \textbf{L = R (verification)} \\
\hline
\endhead

\hline
\endfoot

\hline
\endlastfoot

101 & 38 & 76 & 7676 & 0.0396 & 0.0396 & \text{True} \\
102 & 26 & 2652 & 2652 & 0.0392 & 0.0392 & \text{True} \\
103 & 26 & 5356 & 5356 & 0.0388 & 0.0388 & \text{True} \\
104 & 27 & 1404 & 1404 & 0.0385 & 0.0385 & \text{True} \\
105 & 27 & 1890 & 1890 & 0.0381 & 0.0381 & \text{True} \\
106 & 27 & 2862 & 2862 & 0.0377 & 0.0377 & \text{True} \\
107 & 27 & 5778 & 5778 & 0.0374 & 0.0374 & \text{True} \\
108 & 28 & 1512 & 1512 & 0.037 & 0.037 & \text{True} \\
109 & 28 & 1526 & 3052 & 0.0367 & 0.0367 & \text{True} \\
110 & 28 & 3080 & 3080 & 0.0364 & 0.0364 & \text{True} \\
111 & 28 & 6216 & 6216 & 0.036 & 0.036 & \text{True} \\
112 & 29 & 1624 & 1624 & 0.0357 & 0.0357 & \text{True} \\
113 & 132 & 36 & 22374 & 0.0354 & 0.0354 & \text{True} \\
114 & 29 & 3306 & 3306 & 0.0351 & 0.0351 & \text{True} \\
115 & 29 & 6670 & 6670 & 0.0348 & 0.0348 & \text{True} \\
116 & 30 & 1740 & 1740 & 0.0345 & 0.0345 & \text{True} \\
117 & 30 & 2340 & 2340 & 0.0342 & 0.0342 & \text{True} \\
118 & 30 & 3540 & 3540 & 0.0339 & 0.0339 & \text{True} \\
119 & 30 & 7140 & 7140 & 0.0336 & 0.0336 & \text{True} \\
120 & 31 & 1860 & 1860 & 0.0333 & 0.0333 & \text{True} \\
121 & 33 & 726 & 726 & 0.0331 & 0.0331 & \text{True} \\
122 & 31 & 3782 & 3782 & 0.0328 & 0.0328 & \text{True} \\
123 & 31 & 7626 & 7626 & 0.0325 & 0.0325 & \text{True} \\
124 & 32 & 1984 & 1984 & 0.0323 & 0.0323 & \text{True} \\
125 & 32 & 2400 & 3000 & 0.032 & 0.032 & \text{True} \\
126 & 32 & 4032 & 4032 & 0.0317 & 0.0317 & \text{True} \\
127 & 32 & 8128 & 8128 & 0.0315 & 0.0315 & \text{True} \\
128 & 33 & 2112 & 2112 & 0.0312 & 0.0312 & \text{True} \\
129 & 33 & 2838 & 2838 & 0.031 & 0.031 & \text{True} \\
130 & 33 & 4290 & 4290 & 0.0308 & 0.0308 & \text{True} \\
131 & 33 & 8646 & 8646 & 0.0305 & 0.0305 & \text{True} \\
132 & 34 & 2244 & 2244 & 0.0303 & 0.0303 & \text{True} \\
133 & 34 & 2856 & 3192 & 0.0301 & 0.0301 & \text{True} \\
134 & 34 & 4556 & 4556 & 0.0299 & 0.0299 & \text{True} \\
135 & 34 & 9180 & 9180 & 0.0296 & 0.0296 & \text{True} \\
136 & 35 & 2380 & 2380 & 0.0294 & 0.0294 & \text{True} \\
137 & 35 & 2740 & 3836 & 0.0292 & 0.0292 & \text{True} \\
138 & 35 & 4830 & 4830 & 0.029 & 0.029 & \text{True} \\
139 & 35 & 9730 & 9730 & 0.0288 & 0.0288 & \text{True} \\
140 & 36 & 2520 & 2520 & 0.0286 & 0.0286 & \text{True} \\
141 & 36 & 3384 & 3384 & 0.0284 & 0.0284 & \text{True} \\
142 & 36 & 5112 & 5112 & 0.0282 & 0.0282 & \text{True} \\
143 & 36 & 10296 & 10296 & 0.028 & 0.028 & \text{True} \\
144 & 37 & 2664 & 2664 & 0.0278 & 0.0278 & \text{True} \\
145 & 37 & 3190 & 4070 & 0.0276 & 0.0276 & \text{True} \\
146 & 37 & 5402 & 5402 & 0.0274 & 0.0274 & \text{True} \\
147 & 37 & 10878 & 10878 & 0.0272 & 0.0272 & \text{True} \\
148 & 38 & 2812 & 2812 & 0.027 & 0.027 & \text{True} \\
149 & 70 & 80 & 16688 & 0.0268 & 0.0268 & \text{True} \\
150 & 38 & 5700 & 5700 & 0.0267 & 0.0267 & \text{True} \\
\hline
\end{longtable}


\renewcommand{\thetable}{5} 

\begin{longtable}{|c|c|c|c|c|c|c|}
\caption{Results of Erdős–Straus Conjecture with second formula testing for n from 150 to 200 } \\

\hline
\textbf{n} & \textbf{x} & \textbf{y} & \textbf{z} & \textbf{L} & \textbf{R} & \textbf{L = R (verification)} \\
\hline
\endfirsthead

\hline
\textbf{n} & \textbf{x} & \textbf{y} & \textbf{z} & \textbf{L} & \textbf{R} & \textbf{L = R (verification)} \\
\hline
\endhead

\hline
\endfoot

\hline
\endlastfoot

151 & 38 & 11476 & 11476 & 0.0265 & 0.0265 & \text{True} \\
152 & 39 & 2964 & 2964 & 0.0263 & 0.0263 & \text{True} \\
153 & 39 & 3978 & 3978 & 0.0261 & 0.0261 & \text{True} \\
154 & 39 & 6006 & 6006 & 0.026 & 0.026 & \text{True} \\
155 & 39 & 12090 & 12090 & 0.0258 & 0.0258 & \text{True} \\
156 & 40 & 3120 & 3120 & 0.0256 & 0.0256 & \text{True} \\
157 & 40 & 3768 & 4710 & 0.0255 & 0.0255 & \text{True} \\
158 & 40 & 6320 & 6320 & 0.0253 & 0.0253 & \text{True} \\
159 & 40 & 12720 & 12720 & 0.0252 & 0.0252 & \text{True} \\
160 & 41 & 3280 & 3280 & 0.025 & 0.025 & \text{True} \\
161 & 42 & 1932 & 1932 & 0.0248 & 0.0248 & \text{True} \\
162 & 41 & 6642 & 6642 & 0.0247 & 0.0247 & \text{True} \\
163 & 41 & 13366 & 13366 & 0.0245 & 0.0245 & \text{True} \\
164 & 42 & 3444 & 3444 & 0.0244 & 0.0244 & \text{True} \\
165 & 42 & 4620 & 4620 & 0.0242 & 0.0242 & \text{True} \\
166 & 42 & 6972 & 6972 & 0.0241 & 0.0241 & \text{True} \\
167 & 42 & 14028 & 14028 & 0.024 & 0.024 & \text{True} \\
168 & 43 & 3612 & 3612 & 0.0238 & 0.0238 & \text{True} \\
169 & 44 & 1690 & 2860 & 0.0237 & 0.0237 & \text{True} \\
170 & 43 & 7310 & 7310 & 0.0235 & 0.0235 & \text{True} \\
171 & 43 & 14706 & 14706 & 0.0234 & 0.0234 & \text{True} \\
172 & 44 & 3784 & 3784 & 0.0233 & 0.0233 & \text{True} \\
173 & 45 & 1730 & 3114 & 0.0231 & 0.0231 & \text{True} \\
174 & 44 & 7656 & 7656 & 0.023 & 0.023 & \text{True} \\
175 & 44 & 15400 & 15400 & 0.0229 & 0.0229 & \text{True} \\
176 & 45 & 3960 & 3960 & 0.0227 & 0.0227 & \text{True} \\
177 & 45 & 5310 & 5310 & 0.0226 & 0.0226 & \text{True} \\
178 & 45 & 8010 & 8010 & 0.0225 & 0.0225 & \text{True} \\
179 & 45 & 16110 & 16110 & 0.0223 & 0.0223 & \text{True} \\
180 & 46 & 4140 & 4140 & 0.0222 & 0.0222 & \text{True} \\
181 & 48 & 1086 & 2896 & 0.0221 & 0.0221 & \text{True} \\
182 & 46 & 8372 & 8372 & 0.022 & 0.022 & \text{True} \\
183 & 46 & 16836 & 16836 & 0.0219 & 0.0219 & \text{True} \\
184 & 47 & 4324 & 4324 & 0.0217 & 0.0217 & \text{True} \\
185 & 47 & 5180 & 6580 & 0.0216 & 0.0216 & \text{True} \\
186 & 47 & 8742 & 8742 & 0.0215 & 0.0215 & \text{True} \\
187 & 47 & 17578 & 17578 & 0.0214 & 0.0214 & \text{True} \\
188 & 48 & 4512 & 4512 & 0.0213 & 0.0213 & \text{True} \\
189 & 48 & 6048 & 6048 & 0.0212 & 0.0212 & \text{True} \\
190 & 48 & 9120 & 9120 & 0.0211 & 0.0211 & \text{True} \\
191 & 48 & 18336 & 18336 & 0.0209 & 0.0209 & \text{True} \\
192 & 49 & 4704 & 4704 & 0.0208 & 0.0208 & \text{True} \\
193 & 676 & 52 & 65234 & 0.0207 & 0.0207 & \text{True} \\
194 & 49 & 9506 & 9506 & 0.0206 & 0.0206 & \text{True} \\
195 & 49 & 19110 & 19110 & 0.0205 & 0.0205 & \text{True} \\
196 & 50 & 4900 & 4900 & 0.0204 & 0.0204 & \text{True} \\
197 & 148 & 74 & 29156 & 0.0203 & 0.0203 & \text{True} \\
198 & 50 & 9900 & 9900 & 0.0202 & 0.0202 & \text{True} \\
199 & 50 & 19900 & 19900 & 0.0201 & 0.0201 & \text{True} \\
200 & 51 & 5100 & 5100 & 0.02 & 0.02 & \text{True} \\
\hline
\end{longtable}


\renewcommand{\thetable}{6}

\renewcommand{\arraystretch}{1.15}

\begin{longtable}{|r|r|c|r|r|r|r|}
\caption{Results of Sierpiński Conjecture with second formula testing for $n$ from $2$ to $45$} \\

\hline
\textbf{$n$} & \textbf{$x$} & \textbf{$t$} & \textbf{$y$} & \textbf{$z$} & \textbf{$F$} & \textbf{$q$} \\
\hline
\endfirsthead

\hline
\textbf{$n$} & \textbf{$x$} & \textbf{$t$} & \textbf{$y$} & \textbf{$z$} & \textbf{$F$} & \textbf{$q$} \\
\hline
\endhead

\hline
\endfoot

\hline
\endlastfoot

2  & 1  & $\tfrac12$ & 1   & 2   & $\tfrac14$     & $\tfrac12$ \\
3  & 1  & $2$        & 2   & 6   & $4$            & $2$ \\
4  & 1  & $8$        & 8   & 8   & $0$            & $0$ \\
5  & 2  & $\tfrac45$ & 4   & 4   & $0$            & $0$ \\
6  & 2  & $2$        & 4   & 12  & $16$           & $4$ \\
7  & 5  & $2$        & 2   & 70  & $1156$         & $34$ \\
8  & 2  & $8$        & 16  & 16  & $0$            & $0$ \\
9  & 2  & $36$       & 36  & 36  & $0$            & $0$ \\
10 & 3  & $\tfrac52$ & 10  & 15  & $\tfrac{25}{4}$& $\tfrac52$ \\
11 & 5  & $8$        & 4   & 220 & $11664$        & $108$ \\
12 & 3  & $8$        & 24  & 24  & $0$            & $0$ \\
13 & 3  & $26$       & 26  & 78  & $676$          & $26$ \\
14 & 3  & $84$       & 84  & 84  & $0$            & $0$ \\
15 & 4  & $5$        & 20  & 30  & $25$           & $5$ \\
16 & 4  & $8$        & 32  & 32  & $0$            & $0$ \\
17 & 4  & $17$       & 34  & 68  & $289$          & $17$ \\
18 & 4  & $36$       & 72  & 72  & $0$            & $0$ \\
19 & 4  & $152$      & 152 & 152 & $0$            & $0$ \\
20 & 5  & $8$        & 40  & 40  & $0$            & $0$ \\
21 & 5  & $14$       & 42  & 70  & $196$          & $14$ \\
22 & 5  & $44$       & 44  & 220 & $7744$         & $88$ \\
23 & 5  & $98$       & 70  & 322 & $15876$        & $126$ \\
24 & 5  & $240$      & 240 & 240 & $0$            & $0$ \\
25 & 6  & $12$       & 60  & 60  & $0$            & $0$ \\
26 & 6  & $26$       & 52  & 156 & $2704$         & $52$ \\
27 & 6  & $36$       & 108 & 108 & $0$            & $0$ \\
28 & 6  & $84$       & 168 & 168 & $0$            & $0$ \\
29 & 6  & $348$      & 348 & 348 & $0$            & $0$ \\
30 & 7  & $20$       & 60  & 140 & $1600$         & $40$ \\
31 & 7  & $62$       & 62  & 434 & $34596$        & $186$ \\
32 & 7  & $50$       & 140 & 160 & $100$          & $10$ \\
33 & 7  & $128$      & 176 & 336 & $6400$         & $80$ \\
34 & 7  & $476$      & 476 & 476 & $0$            & $0$ \\
35 & 8  & $35$       & 70  & 280 & $11025$        & $105$ \\
36 & 8  & $36$       & 144 & 144 & $0$            & $0$ \\
37 & 8  & $74$       & 148 & 296 & $5476$         & $74$ \\
38 & 8  & $152$      & 304 & 304 & $0$            & $0$ \\
39 & 8  & $624$      & 624 & 624 & $0$            & $0$ \\
40 & 9  & $30$       & 120 & 180 & $900$          & $30$ \\
41 & 12 & $169$      & 26  & 6396& $10144225$     & $3185$ \\
42 & 9  & $84$       & 252 & 252 & $0$            & $0$ \\
43 & 9  & $258$      & 258 & 774 & $66564$        & $258$ \\
44 & 9  & $792$      & 792 & 792 & $0$            & $0$ \\
45 & 10 & $36$       & 180 & 180 & $0$            & $0$ \\
\hline
\end{longtable}

\renewcommand{\thetable}{7}

\renewcommand{\arraystretch}{1.15}

\begin{longtable}{|r|r|r|r|r|r|r|}
\caption{Results of Sierpiński Conjecture with second formula testing for $n$ from $46$ to $89$} \\

\hline
\textbf{$n$} & \textbf{$x$} & \textbf{$t$} & \textbf{$y$} & \textbf{$z$} & \textbf{$F$} & \textbf{$q$} \\
\hline
\endfirsthead

\hline
\textbf{$n$} & \textbf{$x$} & \textbf{$t$} & \textbf{$y$} & \textbf{$z$} & \textbf{$F$} & \textbf{$q$} \\
\hline
\endhead

\hline
\endfoot

\hline
\endlastfoot

46 & 10 & $98$   & 140  & 644   & $63504$     & $252$ \\
47 & 10 & $188$  & 188  & 940   & $141376$    & $376$ \\
48 & 10 & $240$  & 480  & 480   & $0$         & $0$ \\
49 & 10 & $980$  & 980  & 980   & $0$         & $0$ \\
50 & 11 & $44$   & 220  & 220   & $0$         & $0$ \\
51 & 12 & $17$   & 102  & 204   & $2601$      & $51$ \\
52 & 11 & $128$  & 352  & 416   & $1024$      & $32$ \\
53 & 11 & $512$  & 352  & 1696  & $451584$    & $672$ \\
54 & 11 & $1188$ & 1188 & 1188  & $0$         & $0$ \\
55 & 12 & $55$   & 220  & 330   & $3025$      & $55$ \\
56 & 12 & $84$   & 336  & 336   & $0$         & $0$ \\
57 & 12 & $152$  & 456  & 456   & $0$         & $0$ \\
58 & 12 & $348$  & 696  & 696   & $0$         & $0$ \\
59 & 12 & $1416$ & 1416 & 1416  & $0$         & $0$ \\
60 & 13 & $65$   & 260  & 390   & $4225$      & $65$ \\
61 & 15 & $242$  & 66   & 6710  & $11035684$  & $3322$ \\
62 & 14 & $62$   & 124  & 868   & $138384$    & $372$ \\
63 & 13 & $416$  & 728  & 936   & $10816$     & $104$ \\
64 & 13 & $1664$ & 1664 & 1664  & $0$         & $0$ \\
65 & 14 & $80$   & 280  & 520   & $14400$     & $120$ \\
66 & 14 & $128$  & 352  & 672   & $25600$     & $160$ \\
67 & 15 & $134$  & 134  & 2010  & $879844$    & $938$ \\
68 & 14 & $476$  & 952  & 952   & $0$         & $0$ \\
69 & 14 & $1932$ & 1932 & 1932  & $0$         & $0$ \\
70 & 15 & $84$   & 420  & 420   & $0$         & $0$ \\
71 & 15 & $142$  & 426  & 710   & $20164$     & $142$ \\
72 & 15 & $240$  & 720  & 720   & $0$         & $0$ \\
73 & 15 & $584$  & 876  & 1460  & $85264$     & $292$ \\
74 & 15 & $2220$ & 2220 & 2220  & $0$         & $0$ \\
75 & 16 & $96$   & 480  & 480   & $0$         & $0$ \\
76 & 16 & $152$  & 608  & 608   & $0$         & $0$ \\
77 & 16 & $275$  & 770  & 880   & $3025$      & $55$ \\
78 & 16 & $624$  & 1248 & 1248  & $0$         & $0$ \\
79 & 16 & $2528$ & 2528 & 2528  & $0$         & $0$ \\
80 & 17 & $125$  & 400  & 850   & $50625$     & $225$ \\
81 & 18 & $36$   & 324  & 324   & $0$         & $0$ \\
82 & 24 & $169$  & 52   & 12792 & $40576900$  & $6370$ \\
83 & 20 & $125$  & 100  & 4150  & $4100625$   & $2025$ \\
84 & 17 & $2856$ & 2856 & 2856  & $0$         & $0$ \\
85 & 18 & $240$  & 360  & 2040  & $705600$    & $840$ \\
86 & 18 & $258$  & 516  & 1548  & $266256$    & $516$ \\
87 & 18 & $348$  & 1044 & 1044  & $0$         & $0$ \\
88 & 18 & $792$  & 1584 & 1584  & $0$         & $0$ \\
89 & 18 & $3204$ & 3204 & 3204  & $0$         & $0$ \\

\hline
\end{longtable}

\renewcommand{\thetable}{8} 

\begin{longtable}{|c|c|c|c|c|c|c|}
\caption{Results of the generalized conjecture ($a=6$) with second formula, testing for $n$ from $2$ to $45$} \\

\hline
\textbf{$n$} & \textbf{$x$} & \textbf{$t$} & \textbf{$y$} & \textbf{$z$} & \textbf{$F$} & \textbf{$q$} \\
\hline
\endfirsthead

\hline
\textbf{$n$} & \textbf{$x$} & \textbf{$t$} & \textbf{$y$} & \textbf{$z$} & \textbf{$F$} & \textbf{$q$} \\
\hline
\endhead

\hline
\endfoot

\hline
\endlastfoot

2  & 1  & $1/4$ & 1   & 1    & 0        & 0 \\
3  & 1  & $2/3$ & 2   & 2    & 0        & 0 \\
4  & 1  & 2     & 4   & 4    & 0        & 0 \\
5  & 1  & 10    & 10  & 10   & 0        & 0 \\
6  & 2  & $2/3$ & 4   & 4    & 0        & 0 \\
7  & 3  & 2     & 2   & 42   & 400      & 20 \\
8  & 2  & 2     & 8   & 8    & 0        & 0 \\
9  & 2  & 4     & 12  & 12   & 0        & 0 \\
10 & 2  & 10    & 20  & 20   & 0        & 0 \\
11 & 2  & 44    & 44  & 44   & 0        & 0 \\
12 & 3  & 2     & 12  & 12   & 0        & 0 \\
13 & 3  & 32    & 8   & 312  & 23104    & 152 \\
14 & 3  & 7     & 14  & 42   & 196      & 14 \\
15 & 3  & 10    & 30  & 30   & 0        & 0 \\
16 & 3  & 24    & 48  & 48   & 0        & 0 \\
17 & 3  & 102   & 102 & 102  & 0        & 0 \\
18 & 4  & 4     & 24  & 24   & 0        & 0 \\
19 & 4  & 32    & 16  & 304  & 20736    & 144 \\
20 & 4  & 10    & 40  & 40   & 0        & 0 \\
21 & 4  & 21    & 42  & 84   & 441      & 21 \\
22 & 4  & 44    & 88  & 88   & 0        & 0 \\
23 & 4  & 184   & 184 & 184  & 0        & 0 \\
24 & 5  & 12    & 24  & 120  & 2304     & 48 \\
25 & 5  & 10    & 50  & 50   & 0        & 0 \\
26 & 6  & 32    & 16  & 624  & 92416    & 304 \\
27 & 5  & 30    & 90  & 90   & 0        & 0 \\
28 & 5  & 70    & 140 & 140  & 0        & 0 \\
29 & 5  & 290   & 290 & 290  & 0        & 0 \\
30 & 6  & 10    & 60  & 60   & 0        & 0 \\
31 & 39 & 2     & 6   & 806  & 160000   & 400 \\
32 & 6  & 24    & 96  & 96   & 0        & 0 \\
33 & 6  & 44    & 132 & 132  & 0        & 0 \\
34 & 6  & 102   & 204 & 204  & 0        & 0 \\
35 & 6  & 420   & 420 & 420  & 0        & 0 \\
36 & 7  & 14    & 84  & 84   & 0        & 0 \\
37 & 27 & 32    & 8   & 7992 & 15936064 & 3992 \\
38 & 7  & 76    & 76  & 532  & 51984    & 228 \\
39 & 8  & 54    & 36  & 936  & 202500   & 450 \\
40 & 7  & 140   & 280 & 280  & 0        & 0 \\
41 & 7  & 574   & 574 & 574  & 0        & 0 \\
42 & 8  & 21    & 84  & 168  & 1764     & 42 \\
43 & 8  & 43    & 86  & 344  & 16641    & 129 \\
44 & 8  & 44    & 176 & 176  & 0        & 0 \\
45 & 8  & 80    & 240 & 240  & 0        & 0 \\
\hline
\end{longtable}

\renewcommand{\thetable}{9}

\renewcommand{\thetable}{9}

\begin{longtable}{|c|c|c|c|p{2.5cm}|p{2.5cm}|c|}
\caption{Results of the generalized conjecture ($a=11$) with first formula testing for $n$ from $500$ to $625$} \\

\hline
\textbf{n} & \textbf{x} & \textbf{y} & \textbf{z} & \textbf{L} & \textbf{R} & \textbf{L = R} \\
\hline
\endfirsthead

\hline
\textbf{n} & \textbf{x} & \textbf{y} & \textbf{z} & \textbf{L} & \textbf{R} & \textbf{L = R} \\
\hline
\endhead

\hline
\endfoot

\hline
\endlastfoot

501 & 46 & 11523 & 7682  & 0.021956088 & 0.021956088 & True \\
503 & 46 & 23138 & 11569 & 0.021868787 & 0.021868787 & True \\
504 & 46 & 46368 & 15456 & 0.021825397 & 0.021825397 & True \\
510 & 49 & 1785  & 1666  & 0.021568627 & 0.021568627 & True \\
513 & 50 & 1425  & 1350  & 0.021442495 & 0.021442495 & True \\
514 & 47 & 24158 & 12079 & 0.021400778 & 0.021400778 & True \\
521 & 48 & 8336  & 6252  & 0.021113244 & 0.021113244 & True \\
523 & 48 & 12552 & 8368  & 0.021032505 & 0.021032505 & True \\
525 & 48 & 25200 & 12600 & 0.020952381 & 0.020952381 & True \\
526 & 48 & 50496 & 16832 & 0.020912548 & 0.020912548 & True \\
528 & 51 & 1683  & 1584  & 0.020833333 & 0.020833333 & True \\
531 & 49 & 7434  & 5782  & 0.020715631 & 0.020715631 & True \\
534 & 49 & 13083 & 8722  & 0.020599251 & 0.020599251 & True \\
536 & 49 & 26264 & 13132 & 0.020522388 & 0.020522388 & True \\
537 & 49 & 52626 & 17542 & 0.020484171 & 0.020484171 & True \\
545 & 51 & 3706  & 3270  & 0.020183486 & 0.020183486 & True \\
546 & 50 & 18200 & 10920 & 0.020146520 & 0.020146520 & True \\
547 & 50 & 27350 & 13675 & 0.020109689 & 0.020109689 & True \\
550 & 51 & 5610  & 4675  & 0.020000000 & 0.020000000 & True \\
556 & 51 & 14178 & 9452  & 0.019784173 & 0.019784173 & True \\
558 & 51 & 28458 & 14229 & 0.019713262 & 0.019713262 & True \\
559 & 51 & 57018 & 19006 & 0.019677996 & 0.019677996 & True \\
567 & 52 & 14742 & 9828  & 0.019400353 & 0.019400353 & True \\
569 & 52 & 29588 & 14794 & 0.019332162 & 0.019332162 & True \\
570 & 52 & 59280 & 19760 & 0.019298246 & 0.019298246 & True \\
575 & 54 & 3450  & 3105  & 0.019130435 & 0.019130435 & True \\
576 & 53 & 10176 & 7632  & 0.019097222 & 0.019097222 & True \\
580 & 53 & 30740 & 15370 & 0.018965517 & 0.018965517 & True \\
581 & 54 & 5229  & 4482  & 0.018932874 & 0.018932874 & True \\
584 & 55 & 3212  & 2920  & 0.018835616 & 0.018835616 & True \\
589 & 54 & 15903 & 10602 & 0.018675722 & 0.018675722 & True \\
590 & 54 & 21240 & 12744 & 0.018644068 & 0.018644068 & True \\
591 & 54 & 31914 & 15957 & 0.018612521 & 0.018612521 & True \\
592 & 54 & 63936 & 21312 & 0.018581081 & 0.018581081 & True \\
594 & 55 & 6534  & 5445  & 0.018518519 & 0.018518519 & True \\
596 & 55 & 8195  & 6556  & 0.018456376 & 0.018456376 & True \\
600 & 55 & 16500 & 11000 & 0.018333333 & 0.018333333 & True \\
601 & 56 & 4808  & 4207  & 0.018302829 & 0.018302829 & True \\
602 & 55 & 33110 & 16555 & 0.018272425 & 0.018272425 & True \\
603 & 55 & 66330 & 22110 & 0.018242123 & 0.018242123 & True \\
609 & 56 & 11368 & 8526  & 0.018062397 & 0.018062397 & True \\
610 & 56 & 13664 & 9760  & 0.018032787 & 0.018032787 & True \\
612 & 65 & 780   & 765   & 0.017973856 & 0.017973856 & True \\
613 & 56 & 34328 & 17164 & 0.017944535 & 0.017944535 & True \\
620 & 57 & 11780 & 8835  & 0.017741935 & 0.017741935 & True \\
622 & 57 & 17727 & 11818 & 0.017684887 & 0.017684887 & True \\
624 & 57 & 35568 & 17784 & 0.017628205 & 0.017628205 & True \\

\hline
\end{longtable}

\appendix

\section*{Appendix: Python Code for Testing the Generalized Erdős-Straus Conjecture (First Formula)}

The following code allows testing the generalized Erdős-Straus conjecture numerically for different values of \(a\) and \(n\). It provides direct output of the results without generating a CSV file.

\lstset{
    language=Python,
    basicstyle=\ttfamily\footnotesize,
    keywordstyle=\color{blue},
    commentstyle=\color{gray},
    stringstyle=\color{red},
    breaklines=true,
    numbers=left,
    numberstyle=\tiny,
    frame=single,
    showstringspaces=false
}

\begin{lstlisting}
# Python implementation of the first formula for the Generalized Erdos-Straus Conjecture


def generalized_erdos_straus_conjecture(a, n_min, n_max):
    """
    Tests the generalized Erdos-Straus conjecture for a >= 4.
.
    This function prints results directly instead of saving them to a file.
    
    :param a: Integer, the numerator in the equation a/n = 1/x + 1/y + 1/z.
    :param n_min: Integer, the starting value of n.
    :param n_max: Integer, the ending value of n.
    """
    epsilon = 1e-9  # Tolerance for floating-point comparison
    for n in range(n_min, n_max + 1):
        solution_found = False
        L = a / n  # Left-hand side of the equation
        for x in range(n // a + 1, 10 * n):  # Searching for x > n/a
            lambda_value = a * x - n - 1
            if lambda_value <= 0:
                continue
            y = (2 * n * x) // lambda_value
            z = (2 * n * x) // (lambda_value + 2)

            # Ensure y and z are integers
            if (2 * n * x) % lambda_value == 0 and (2 * n * x) % (lambda_value + 2) == 0:
                R = 1/x + 1/y + 1/z  # Right-hand side of the equation
                is_true = abs(L - R) < epsilon  # Floating-point comparison with tolerance
                print(f"For n = {n}, x = {x}, y = {y}, z = {z}: L = {L:.4f}, R = {R:.4f}, Conjecture holds: {is_true}")
                solution_found = True
                break

# Example usage:
a = 4 #(for the classic Erdos_Strauss Conjecture)
n_min = 2
n_max = 100
generalized_erdos_straus_conjecture(a, n_min, n_max)
\end{lstlisting}

\newpage

\section*{Appendix: Python Code for Testing the Generalized Erdős-Straus Conjecture (Second Formula)}

The following code tests the generalized Erdős-Straus conjecture using the second formula independently. It prints results directly based on the existence of a perfect square.

\begin{lstlisting}
import math
import time
from fractions import Fraction

print("VERIFICATION OF OUR CONJECTURE -- CASE a = 4")
print("="*70)

def verify_our_conjecture(n_min=2, n_max=500,
                     max_x_factor=300,
                     max_t_range=500,
                     max_q=10):

    captured = 0
    failed_n_list = []

    start_time = time.time()
    a = 4# You can change de value of a to test another conjecture

    for n in range(n_min, n_max + 1):
        found = False

        x_start = n // a + 1
        x_limit = max_x_factor * n

        for x in range(x_start, x_limit + 1):
            d = a * x - n
            if d == 0:
                continue

            denom = d * d
            t_min = max(1, (2 * n * x) // denom)

            # =========================
            # PHASE 1 : integer t
            # =========================
            for t in range(t_min, t_min + max_t_range + 1):
                F = t*t*d*d - 2*t*x*n
                if F < 0:
                    continue

                q = math.isqrt(F)
                if q*q == F:
                    y = t*d - q
                    z = t*d + q
                    if y > 0 and z > 0 and y == int(y) and z == int(z):
                        if a*x*y*z == n*(y*z + x*z + x*y):
                            captured += 1
                            found = True
                            print(f"n={n:4d} | x={x}, t={t} [INTEGER], y={y}, z={z}, F={F}, q={q}")
                            break
            if found:
                break

        # =========================
        # PHASE 2 : rational t
        # =========================
        if not found:
            for x in range(x_start, x_limit + 1):
                d = a * x - n
                if d == 0:
                    continue

                for q_den in range(2, max_q+1):
                    for p in range(1, max_q*30):
                        t = Fraction(p, q_den)
                        F = t*t*d*d - 2*t*x*n
                        if F < 0:
                            continue

                        F_num = F.numerator
                        F_den = F.denominator
                        if math.isqrt(F_num)**2 == F_num and math.isqrt(F_den)**2 == F_den:
                            qrat = Fraction(math.isqrt(F_num), math.isqrt(F_den))
                            y = t*d - qrat
                            z = t*d + qrat
                            # === require y and z to be integers ===
                            if y.denominator == 1 and z.denominator == 1 and y > 0 and z > 0:
                                if a*x*int(y)*int(z) == n*(int(y)*int(z) + x*int(z) + x*int(y)):
                                    captured += 1
                                    found = True
                                    print(f"n={n:4d} | x={x}, t={t} [RATIONAL], y={y}, z={z}, F={F}, q={qrat}")
                                    break
                    if found:
                        break
                if found:
                    break

        if not found:
            print(f"n={n:4d} : X FAILED")

            failed_n_list.append(n)

    elapsed_time = time.time() - start_time
    total = n_max - n_min + 1
    success_rate = 100 * captured / total

    print("\n" + "="*70)
    print("FINAL REPORT a = 4")
    print("="*70)
   print(f"Tested interval : n in [{n_min}, {n_max}]")

    print(f"Success : {captured}/{total} ({success_rate:.2f}%)")
    print(f"Time : {elapsed_time:.2f}s")

    if failed_n_list:
        print("Detected failures :", failed_n_list)
    else:
        print(" COMPLETE SUCCESS")

    return    



# DIRECT EXECUTION


verify_our_conjecture(n_min=2, n_max=500)#You can vary the check interval for a fixed a
\end{lstlisting}

\newpage

\section*{Appendix: Python Code for Testing the Classic Erdos-Straus Conjecture for the Values Not Proven by Mordell (Second Formula)}

The following code tests the generalized Erdős-Straus conjecture using the second formula specifically for values of \(n\) where Mordell's method was inconclusive. 

\begin{lstlisting}

import math
import time

print("VERIFICATION OF ERDOS-STRAUS CONJECTURE - a = 4 (Our Method)")

print("Testing n in quadratic residues modulo 840")
print("=" * 70)

# Parameters
a = 4
MOD = 840
RESIDUES = {(r * r) % MOD for r in [1, 11, 13, 17, 19, 23]}

def generate_admissible_n(n_min, n_max):
    """Generate list of admissible n (quadratic residues mod 840) in [n_min, n_max]"""
    admissible_list = []
    for n in range(n_min, n_max + 1):
        if (n % MOD) in RESIDUES:
            admissible_list.append(n)
    return admissible_list

def verify_with_our_method_for_list(n_list, max_x_factor=300, max_t_range=500):
    """
    Verify the conjecture using our method for a given list of n.
    Returns: solutions dict, failures list, elapsed time
    """

    captured = 0
    failed_n_list = []
    solutions = {}

    start_time = time.time()
    total_n = len(n_list)

    print(f"\nTesting {total_n} admissible n")
    print(f"Parameters: max_x_factor={max_x_factor}, max_t_range={max_t_range}")
    print("-" * 70)

    for idx, n in enumerate(n_list, 1):
        found = False

        # Search for x
        x_start = n // a + 1
        x_limit = max_x_factor * n

        for x in range(x_start, x_limit + 1):
            if found:
                break

            d = a * x - n
            if d == 0:
                continue

            # Compute t_min to ensure F >= 0
            denom = d * d
            t_min = max(1, (2 * n * x) // denom)

            # Search for t
            for t in range(t_min, t_min + max_t_range + 1):
                F = t * t * d * d - 2 * t * x * n

                if F < 0:
                    continue

                q = math.isqrt(F)
                if q * q == F:
                    y = t * d - q
                    z = t * d + q

                    if y > 0 and z > 0:
                        left = a * x * y * z
                        right = n * (y*z + x*z + x*y)

                        if left == right:
                            solutions[n] = (x, t, y, z)
                            captured += 1
                            found = True

                            # Print solution
                            print(f"n={n:4d} : 0k x={x:5d}, t={t:5d}, y={y:6d}, z={z:6d}")
                            break

        if not found:
            failed_n_list.append(n)
            print(f"n={n:4d} : FAILED")


    elapsed_time = time.time() - start_time

    return solutions, failed_n_list, elapsed_time

def intensive_search(failed_list, max_x_factor=1000, max_t_range=2000):
    """Intensive search for failed n"""

    print(f"\n" + "="*70)
    print(f"INTENSIVE SEARCH FOR {len(failed_list)} FAILED N")
    print("="*70)

    a = 4
    new_solutions = {}

    for idx, n in enumerate(failed_list, 1):
        print(f"Failed {idx}/{len(failed_list)} - n={n} : ", end="", flush=True)

        found = False
        x_start = n // a + 1
        x_limit = max_x_factor * n

        for x in range(x_start, min(x_limit, 50000) + 1):
            if found:
                break

            d = a * x - n
            if d <= 0:
                continue

            t_min = max(1, (2 * n * x) // (d * d))
            t_max = t_min + max_t_range

            for t in range(t_min, min(t_max, t_min + 5000) + 1):
                F = t * t * d * d - 2 * t * x * n

                if F < 0:
                    continue

                q = math.isqrt(F)
                if q * q == F:
                    y = t * d - q
                    z = t * d + q

                    if y > 0 and z > 0:
                        left = a * x * y * z
                        right = n * (y*z + x*z + x*y)

                        if left == right:
                            new_solutions[n] = (x, t, y, z)
                            found = True
                            print(f"ok (x={x}, t={t}, y={y}, z={z})")
                            break

        if not found:
            print("FAILED")


    return new_solutions

# ========================================
# MAIN EXECUTION
# ========================================
if __name__ == "__main__":
    # MODIFY THESE VALUES TO TEST DIFFERENT INTERVALS
    n_min = 5000      # YOU CAN CHANGE THIS
    n_max = 10000   # YOU CAN CHANGE THIS

    # Step 1: Generate admissible n in the interval
    print(f"\nGenerating admissible n in [{n_min}, {n_max}]")
    print("-" * 70)

    admissible_list = generate_admissible_n(n_min, n_max)

    if len(admissible_list) == 0:
        print("No admissible n in this interval!")
        exit()

    print(f"Total admissible n: {len(admissible_list)}")
    print(f"\nList of admissible n ({len(admissible_list)} numbers):")
    # Print admissible n in groups of 20 for readability
    for i in range(0, len(admissible_list), 20):
        print(f"  {admissible_list[i:i+20]}")

    # Step 2: Verify using our method
    solutions, failed, elapsed = verify_with_our_method_for_list(
        n_list=admissible_list,
        max_x_factor=300,
        max_t_range=500
    )

    # Step 3: Final report
    print(f"\n" + "="*70)
    print("FINAL REPORT")
    print("="*70)

    success_rate = 100 * len(solutions) / len(admissible_list)

    print(f"Interval tested        : [{n_min}, {n_max}]")
    print(f"Admissible n tested    : {len(admissible_list)}")
    print(f"Execution time         : {elapsed:.2f} seconds")
    print(f"Solutions found        : {len(solutions)}/{len(admissible_list)} ({success_rate:.2f}%)")
    print(f"Failures               : {len(failed)}")

    if failed:
        print(f"\nList of failures ({len(failed)} numbers):")
        for i in range(0, len(failed), 20):
            print(f"  {failed[i:i+20]}")

        # Step 4: Intensive search for failures
        new_solutions = intensive_search(
            failed_list=failed,
            max_x_factor=1000,
            max_t_range=2000
        )

        if new_solutions:
            print(f"\n" + "="*70)
            print(f"INTENSIVE SEARCH RESULTS")
            print("="*70)

            print(f"New solutions found: {len(new_solutions)}")
            for n in sorted(new_solutions.keys()):
                x, t, y, z = new_solutions[n]
                print(f"  n={n:4d} : x={x:6d}, t={t:6d}, y={y:8d}, z={z:8d}")

            # Update statistics
            solutions.update(new_solutions)
            remaining_failures = [n for n in failed if n not in new_solutions]
            new_total = len(solutions)
            new_rate = 100 * new_total / len(admissible_list)

            print(f"\nUpdated totals:")
            print(f"  Total solutions : {new_total}/{len(admissible_list)} ({new_rate:.2f}%)")
            print(f"  Remaining failures : {len(remaining_failures)}")

            if remaining_failures:
                print(f"\nRemaining failures ({len(remaining_failures)} numbers):")
                for i in range(0, len(remaining_failures), 20):
                    print(f"  {remaining_failures[i:i+20]}")
        else:
            print(f"\nNo new solutions found despite intensive search")
    else:
        print(f"\n COMPLETE SUCCESS! All admissible n have a solution.")

    print(f"\n" + "="*70)
    print("END OF PROGRAM")
    print("="*70)

    # Return the admissible list as well (useful for further analysis)
    admissible_list_result = admissible_list

\end{lstlisting}

Users can execute this script in Python to test numerical values and verify the conjecture for different cases.

\end{document}